\documentclass[reqno]{amsart}
\usepackage[top=.75in, bottom=.75in, left=1in, right=1in]{geometry}
\usepackage[dvipsnames]{xcolor}
\usepackage{mathdots,graphicx,mathrsfs,enumerate,wrapfig,verbatim,placeins,latexsym,float, amsmath,longtable, amssymb,epsfig,breqn, array, multirow, mathtools,thmtools,thm-restate,multicol,diagbox,nicematrix,url,bbold,color,subcaption,hyperref} 
\hypersetup{
    colorlinks = true,
    citecolor = ForestGreen,
    linkcolor = blue,
}
\usepackage{pgf,tikz,pgfplots,tikz-cd}
\pgfplotsset{compat=1.15}
\usepackage[export]{adjustbox}

\usepackage[style=alphabetic]{biblatex}
\definecolor{lightgray}{gray}{0.85}

\newcommand{\C}{\mathbb{C}}

\newcommand{\F}{\mathcal{F}}

\renewcommand{\t}{\tau}

\newcommand{\x}{\mathbf{x}}

\newcommand{\y}{\mathbf{y}}
\newcommand{\Z}{\mathbb{Z}}
\newcommand{\csm}{\textrm{csm}}

\newcommand{\perm}{\textrm{perm}}

\newcommand{\rank}{\textrm{rank}}

\def\Hom{\textrm{Hom}}

\newcommand\plaqctr[1]{\begin{tikzpicture}[baseline={([yshift=-\the\dimexpr\fontdimen22\textfont2\relax]current  bounding  box.center)}]\plaq{#1}\end{tikzpicture}}
\makeatletter
\newcommand{\gettikzxy}[3]{
  \tikz@scan@one@point\pgfutil@firstofone#1\relax
\pgfmathsetmacro{#2}{\the\pgf@x/\linkpatternunit}
\pgfmathsetmacro{#3}{\the\pgf@y/\linkpatternunit}
}
\tikzset{label anchor/.code={%
    \let\tikz@auto@anchor=\pgfutil@empty
    \def\tikz@anchor{#1}
  },
  label anchor/.default=center
}
\makeatother
\tikzset{arrow/.style={postaction={decorate,thick,decoration={markings,mark = at position #1 with {\arrow{>}}}}},arrow/.default=0.5}
\tikzset{invarrow/.style={postaction={decorate,thick,decoration={markings,mark = at position #1 with {\arrow{<}}}}},invarrow/.default=0.5}
\newdimen\linkpatternunit%
\newcount\linkpatternsize%
\newcount\lpsize
\newif\iflinkpatterninverted
\newif\iflinkpatterntikzstarted
\newif\iflinkpatternboxed
\newif\iflinkpatternaxis
\newif\iflinkpatternstraightlines
\newif\iflinkpatternnumbered
\newif\iflinkpatternalias
\newif\iflinkpatternnode
\newif\iflinkpatterncentered
\newcount\linkpatternfused
\pgfkeys{/linkpattern/.cd,centered/.is if=linkpatterncentered,inverted/.is if=linkpatterninverted,numbered/.is if=linkpatternnumbered,tikzstarted/.is if=linkpatterntikzstarted,straight lines/.is if=linkpatternstraightlines,boxed/.is if=linkpatternboxed,axis/.is if=linkpatternaxis,vertexcolor/.store in=\linkpatternvertexcolor,
edgecolor/.store in=\linkpatternedgecolor,
boxcolor/.code={\linkpatternboxedtrue\def\linkpatternboxcolor{#1}},tikzoptions/.style={every linkpattern/.append style={#1}},
size/.code={\linkpatternsize=#1},
numbering0/.code={\def\lpnumbering{#1}\def\linkpatternnumbering{#1}},
numbering/.style={numbered,numbering0={#1}},
unit/.code={\linkpatternunit=#1},
height/.store in=\linkpatternheight,
shape/.store in=\linkpatternshape,
looseness/.store in=\linkpatternlooseness,
squareness/.store in=\linkpatternsquareness,
extra space/.store in=\linkpatternextraspace,
width/.store in=\linkpatternwidth,
alias/.is if=linkpatternalias,
pos/.store in=\linkpatternpos,
labeloptions/.style={labeloptionslist/.append style={#1}},
labeloptionslist/.style={inner sep=2pt,font=\scriptsize,
execute at begin node=$,execute at end node=$,
label anchor={#1+180}},
nodeon/.is if=linkpatternnode,
node/.style={nodeon,nodeoptionslist/.append style={#1}},
nodeoptionslist/.style={},
pipedream/.style={shape=pipedream,looseness=0,straight lines,numbering0=tangle},
tangle/.style={shape=tangle,numbering0=tangle},
every linkpattern/.style={x=\linkpatternunit,y=\linkpatternunit},
fused/.code={\linkpatternfused=#1},
vertex/.style={circle,thin,solid,draw=black,fill=\linkpatternvertexcolor,inner sep=1.5pt,draw opacity=1,transform shape},
edge/.style={very thick,solid,draw=\linkpatternedgecolor,draw opacity=1}}
\linkpatterncenteredfalse%
\linkpatterninvertedfalse%
\linkpatternnumberedfalse%
\linkpatterntikzstartedfalse%
\linkpatternboxedfalse%
\linkpatternaxistrue%
\linkpatternaliastrue%
\linkpatternstraightlinesfalse%
\def\linkpatternlooseness{0.2}
\def\linkpatternsquareness{0.35}
\def\linkpatternvertexcolor{red}%
\def\linkpatternedgecolor{blue}%
\def\linkpatternboxcolor{none}%
\def\linkpatternheight{0}
\def\linkpatternwidth{0}
\def\linkpatternshape{default}
\def\linkpatternnumbering{default}
\def\linkpatternpos{(0,0)}
\def\linkpatternextraspace{0}
\linkpatternnodefalse
\def\firstchar#1#2\empty{#1}%
\def\linkpatterndo#1#2{
\edef\param{\csname linkpattern#2\endcsname}
\edef\firstcharparam{\expandafter\firstchar\param\empty}
\expandafter\ifcat\firstcharparam a
\expandafter\ifx\csname linkpattern#1\param\endcsname\relax
\csname linkpattern#1unknown\endcsname
\else
\csname linkpattern#1\csname linkpattern#2\endcsname\endcsname
\fi
\else
\csname linkpattern#1unknown\endcsname
\fi
}%
\def\linkpatterncoordtangle{\ifnum\x>\lphalfsize\pgfmathparse{\lpsize+1-\x}\xdef\lpcoordx{\pgfmathresult}\xdef\lpcoordy{\lpheight}\xdef\lpangle{270}\else\xdef\lpcoordx{\x}\xdef\lpcoordy{-\lpheight}\xdef\lpangle{90}\fi}
\def\linkpatterncoordpipedream{\ifnum\x>\lphalfsize\pgfmathparse{\lpsize+1-\x-0.5}\xdef\lpcoordx{\pgfmathresult}\xdef\lpcoordy{0}\xdef\lpangle{270}\else\pgfmathparse{0.5-\x}\xdef\lpcoordy{\pgfmathresult}\xdef\lpcoordx{0}\xdef\lpangle{0}\fi}
\def\linkpatterncoordrectangle{
\ifnum\x>\lptqsize
\pgfmathparse{\lpsize+1-\x-0.5}\xdef\lpcoordx{\pgfmathresult}\xdef\lpcoordy{0}\xdef\lpangle{270}
\else\ifnum\x>\lphalfsize
\pgfmathparse{\x-\lptqsize-0.5}\xdef\lpcoordy{\pgfmathresult}\xdef\lpcoordx{\linkpatternwidth}\xdef\lpangle{180}
\else\ifnum\x>\linkpatternheight
\pgfmathparse{\x-\linkpatternheight-0.5}\xdef\lpcoordx{\pgfmathresult}\xdef\lpcoordy{-\linkpatternheight}\xdef\lpangle{90}
\else
\pgfmathparse{0.5-\x}\xdef\lpcoordy{\pgfmathresult}\xdef\lpcoordx{0}\xdef\lpangle{0}
\fi\fi\fi
}%
\def\linkpatternsetsizeunknown{
\global\lpsize=\linkpatternsize
\if\linkpatternheight0
\xdef\maxsep{0}
\foreach \x/\xx in \mylist%
{%
\edef\tempx{\withoutprime{\x}}
\edef\tempxx{\withoutprime{\xx}}
\pgfmathparse{max(\maxsep,abs(\tempx-\tempxx))}
\xdef\maxsep{\pgfmathresult}
}%
\pgfmathparse{0.25+0.8*\linkpatternsquareness*\maxsep}
\xdef\lpheight{\pgfmathresult}
\else
\xdef\lpheight{\linkpatternheight}
\fi
}
\newcount\tempsize
\def\linkpatternrightmostunknown{
\global\lpsize=0
\global\tempsize=0
\foreach\x/\labx in \linkpatternnumbering
{
\edef\tempx{\withoutprime{\x}}
\ifnum\lpsize<\tempx\global\lpsize=\tempx\fi
\global\advance\tempsize by 1
}
\ifnum\tempsize>\lpsize\global\lpsize=\tempsize\fi
}%
\def\linkpatternrightmostdefault{
\global\lpsize=0
\global\tempsize=0
\foreach \x/\y in \mylist
{
\edef\tempx{\withoutprime{\x}}
\ifnum\lpsize<\tempx\global\lpsize=\tempx\fi
\ifx\x\y
\global\advance\tempsize by 1
\else
\edef\tempy{\withoutprime{\y}}
\ifnum\lpsize<\tempy\global\lpsize=\tempy\fi%
\global\advance\tempsize by 2
\fi
}
\ifnum\tempsize>\lpsize\global\lpsize=\tempsize\fi
}%
\def\linkpatternrightmosttangle{
\global\lpsize=0
\global\tempsize=0
\foreach \x/\y in \mylist
{
\edef\tempx{\withoutprime{\x}}
\ifnum\lpsize<\tempx\global\lpsize=\tempx\fi
\ifx\x\y
\global\advance\tempsize by 1
\else
\edef\tempy{\withoutprime{\y}}
\ifnum\lpsize<\tempy\global\lpsize=\tempy\fi%
\global\advance\tempsize by 2
\fi
}
\global\advance\lpsize by\lpsize
\ifnum\tempsize>\lpsize\global\lpsize=\tempsize\fi
}%

\newcommand\linkpattern[2][]{
{
\pgfkeys{/linkpattern/.cd,#1}
\edef\mylist{#2}
\def\primetest##1'{}%
\def\hasaprime##1{\expandafter\primetest##1''}
\def\internalwithoutprime##1'{##1}%
\def\withoutprime##1{\if\hasaprime##1 %
\expandafter\internalwithoutprime##1\else ##1\fi}%
\iflinkpatternnumbered%
\iflinkpatterninverted
\tikzset{/linkpattern/lbl/.style n args={3}{label={[/linkpattern/labeloptionslist=-##1,##3] ##1:##2}}}%
\else%
\tikzset{/linkpattern/lbl/.style n args={3}{label={[/linkpattern/labeloptionslist=##1,##3] ##1:##2}}}%
\fi%
\else%
\tikzset{/linkpattern/lbl/.style={}}%
\fi%
\tikzifinpicture{\linkpatterntikzstartedtrue%
\begin{scope}[shift=\linkpatternpos,/linkpattern/every linkpattern]
}{%
\linkpatterntikzstartedfalse%
\iflinkpatterncentered
\begin{tikzpicture}[baseline=(current  bounding  box.center),/linkpattern/every linkpattern]%
\else
\begin{tikzpicture}[baseline=0,/linkpattern/every linkpattern]%
\fi
}%
\begin{scope}[local bounding box=link pattern box]
\iflinkpatterninverted%
\begin{scope}[yscale=-1]%
\fi%
\linkpatterndo{setsize}{shape}
\ifnum\lpsize=0
\linkpatterndo{rightmost}{numbering}
\fi
\pgfmathtruncatemacro{\lphalfsize}{\lpsize/2}
\linkpatterndo{numbering}{numbering}
\iflinkpatternboxed
\linkpatterndo{drawbox}{shape}
\else
\iflinkpatternaxis
\linkpatterndo{drawaxis}{shape}
\fi
\fi
\foreach\xx/\xlab/\opt in \lpnumbering
{
\ifx\xlab\opt\def\opt{}\fi
\if\hasaprime\xx %
\pgfmathtruncatemacro{\xx}{\lpsize+1-\withoutprime{\xx}}
\fi
\ifnum\linkpatternfused>1
\pgfmathsetmacro{\x}{0.4*(0.5+\linkpatternfused*(0.5+floor((\xx-1)/\linkpatternfused)))+0.6*\xx}
\else
\def\x{\xx}
\fi
\linkpatterndo{coord}{shape}
\iflinkpatternalias\def\xlabb{\xlab}\else\def\xlabb{\xx}\fi
\path (\lpcoordx,\lpcoordy) coordinate[/linkpattern/vertex,/linkpattern/lbl={\lpangle+180}{\xlab}{\opt},alias=v\xlabb] (v\xx) ++(\lpangle:\linkpatternunit) coordinate[alias=vv\xlabb] (vv\xx); 
}
\foreach \a/\b/\c in \mylist
{
\if\hasaprime\a %
\pgfmathtruncatemacro{\a}{\lpsize+1-\withoutprime{\a}}
\fi
\ifx\b\c\def\c{}\fi
\draw[/linkpattern/edge]
\ifx\a\b
(v\a)
\c
--
++(0,\lpheight);
\else
\pgfextra{
\if\hasaprime\b %
\pgfmathtruncatemacro{\b}{\lpsize+1-\withoutprime{\b}}
\fi
\gettikzxy{(v\a)}{\ax}{\ay}
\gettikzxy{(v\b)}{\bx}{\by}
\gettikzxy{(vv\a)}{\axx}{\ayy}
\gettikzxy{(vv\b)}{\bxx}{\byy}
\pgfmathsetmacro{\dist}{sqrt((\ax-\bx)*(\ax-\bx)+(\ay-\by)*(\ay-\by))}
\pgfmathsetmacro{\abx}{(\axx-\ax)*\dist*\linkpatternsquareness+(\bx-\ax)*\linkpatternlooseness)}
\pgfmathsetmacro{\aby}{(\ayy-\ay)*\dist*\linkpatternsquareness+(\by-\ay)*\linkpatternlooseness)}
\pgfmathsetmacro{\bax}{(\bxx-\bx)*\dist*\linkpatternsquareness+(\ax-\bx)*\linkpatternlooseness)}
\pgfmathsetmacro{\bay}{(\byy-\by)*\dist*\linkpatternsquareness+(\ay-\by)*\linkpatternlooseness)}
}
(v\a)
\c
\iflinkpatternstraightlines
\pgfextra{
\pgfmathsetmacro{\t}{((\ax-\bx)*\bay-(\ay-\by)*\bax)/(\aby*\bax-\abx*\bay)}
\pgfmathsetmacro{\abx}{\t*\abx}
\pgfmathsetmacro{\aby}{\t*\aby}
}
[rounded corners=0.2\linkpatternunit] -- ++(\abx,\aby) -- (v\b);
\else
.. controls ++(\abx,\aby) and ++(\bax,\bay) .. 
\fi
(v\b);
\fi
}
\end{scope}
\iflinkpatternnode
\node[fit=(link pattern box),/linkpattern/nodeoptionslist] {};
\fi
\iflinkpatterninverted
\end{scope}
\fi
\iflinkpatterntikzstarted
\end{scope}
\else%
\end{tikzpicture}%
\fi%
}}%
\newcommand\tanglelinkpattern[3][]{%
{
\pgfkeys{/linkpattern/.cd,#1}
\iflinkpatterninverted
\begin{tikzpicture}[/linkpattern/every linkpattern,baseline=\linkpatternunit]%
\else
\begin{tikzpicture}[/linkpattern/every linkpattern,baseline=-\linkpatternunit]%
\fi
\linkpattern[#1,tikzstarted,numbered=false]{#3}
\pgfmathtruncatemacro{\lptempsize}{2*\linkpatternsize}
\iflinkpatterninverted
\begin{scope}[yshift=0.5*\linkpatternunit]
\else
\begin{scope}[yshift=-0.5*\linkpatternunit]
\fi
\linkpattern[tangle,#1,tikzstarted,size=\lptempsize,
numbering=halftangle,
height=0.5]{#2}
\end{scope}
\end{tikzpicture}%
}}
\newcommand\diag[4][]{%
\pgfkeys{/linkpattern/.cd,#1}
\iflinkpatterntikzstarted\else%
\begin{tikzpicture}[scale=0.5]
\fi%
\iflinkpatterninverted%
\begin{scope}[yscale=-1]%
\fi%
\draw (0,0) grid (#2,#3);
\edef\mylist{#4}
\foreach\y/\x/\z in \mylist
{
\ifx\x\z
\draw[decorate,decoration={zigzag,
amplitude=1pt,segment length=5pt}]
(\x-0.5,#3) -- (\x-0.5,\y-0.5) node[circle,fill=black,inner sep=2pt] {} -- (#2,\y-0.5);
\else
\node at (\x-0.5,\y-0.5) {$\z$};
\fi
}
\iflinkpatterninverted
\end{scope}
\fi
\iflinkpatterntikzstarted\else%
\end{tikzpicture}%
\fi%
}
\makeatletter
\tikzset{circle split part fill/.style  args={#1,#2}{%
 alias=tmp@name,
  postaction={%
    insert path={
     \pgfextra{%
     \pgfpointdiff{\pgfpointanchor{\pgf@node@name}{center}}%
                  {\pgfpointanchor{\pgf@node@name}{east}}%
     \pgfmathsetmacro\insiderad{\pgf@x}
      \fill[#1] (\pgf@node@name.base) ([xshift=-\pgflinewidth]\pgf@node@name.east) arc
                          (0:180:\insiderad-\pgflinewidth)--cycle;
      \fill[#2] (\pgf@node@name.base) ([xshift=\pgflinewidth]\pgf@node@name.west)  arc
                           (180:360:\insiderad-\pgflinewidth)--cycle;                    }}}}}  
 \makeatother
\tikzset{bdot/.style={circle,circle split,draw,circle split part fill={black,white},thin,inner sep=1pt}}%
\tikzset{wdot/.style={circle,circle split,draw,circle split part fill={white,black},thin,inner sep=1pt}}%
\newdimen{\loopcellsize}
\tikzset{bgplaq/.style={draw=black,fill=\linkpatternboxcolor}}
\def\plaqwest{}
\def\plaqeast{}
\def\plaqnorth{}
\def\plaqsouth{}
\pgfkeys{/linkpattern/.cd,west/.store in=\plaqwest,east/.store in=\plaqeast,north/.store in=\plaqnorth,south/.store in=\plaqsouth}
\def\plaqname{plaq}

\newcommand\plaq[2][]{
\node[bgplaq,rectangle,draw,minimum size=\loopcellsize,transform shape] (\plaqname) {};
\useasboundingbox;
\pgfkeys{/linkpattern/.cd,#1}
\ifx#2\empty\else
\begin{scope}[x=\loopcellsize,y=\loopcellsize]
\csname plaq#2\endcsname
\end{scope}\fi
}
\newcommand\nbplaq[2][]{
\node[bgplaq,rectangle,draw=none,fill=none,minimum size=\loopcellsize,transform shape] (\plaqname) {};
\useasboundingbox;
\pgfkeys{/linkpattern/.cd,#1}
\ifx#2\empty\else
\begin{scope}[x=\loopcellsize,y=\loopcellsize]
\csname plaq#2\endcsname
\end{scope}\fi
}
\newcommand\ecplaq[3][]{
\node[bgplaq,rectangle,draw,minimum size=\loopcellsize,transform shape] (\plaqname) {};
\useasboundingbox;
\pgfkeys{/linkpattern/.cd,#1}
\begin{scope}[x=\loopcellsize,y=\loopcellsize]
\def\linkpatternedgecolor{#3}
\csname plaq#2\endcsname
\end{scope}
}
\newcommand\wecplaq[3][]{
\node[bgplaq,rectangle,draw,fill=none,minimum size=\loopcellsize,transform shape] (\plaqname) {};
\useasboundingbox;
\pgfkeys{/linkpattern/.cd,#1}
\begin{scope}[x=\loopcellsize,y=\loopcellsize]
\def\linkpatternedgecolor{#3}
\csname plaq#2\endcsname
\end{scope}
}
\newcommand\ecbcplaq[4][]{
\node[bgplaq,rectangle,draw,fill=#4,minimum size=\loopcellsize,transform shape] (\plaqname) {};
\useasboundingbox;
\pgfkeys{/linkpattern/.cd,#1}
\begin{scope}[x=\loopcellsize,y=\loopcellsize]
\def\linkpatternedgecolor{#3}
\csname plaq#2\endcsname
\end{scope}
}

\newcommand\wplaq[2][]{
\node[bgplaq,rectangle,draw,fill=white,minimum size=\loopcellsize,transform shape] (\plaqname) {};
\useasboundingbox;
\pgfkeys{/linkpattern/.cd,#1}
\ifx#2\empty\else
\begin{scope}[x=\loopcellsize,y=\loopcellsize]
\csname plaq#2\endcsname
\end{scope}\fi
}

\newcommand\gplaq[2][]{
\node[bgplaq,rectangle,draw,fill=lightgray,minimum size=\loopcellsize,transform shape] (\plaqname) {};
\useasboundingbox;
\pgfkeys{/linkpattern/.cd,#1}
\ifx#2\empty\else
\begin{scope}[x=\loopcellsize,y=\loopcellsize]
\csname plaq#2\endcsname
\end{scope}\fi
}

\tikzset{loop/.code={\def\plaqname{loop-\the\pgfmatrixcurrentrow-\the\pgfmatrixcurrentcolumn}},loop/.append style={matrix,row sep={\loopcellsize,between origins},column sep={\loopcellsize,between origins}}}
\def\linkpatternboxcolor{pink!50!white}

\numberwithin{equation}{section}

\newtheorem{theorem}{Theorem}[section]
\newtheorem*{theorem*}{Theorem}
\newtheorem{lemma}[theorem]{Lemma}
\newtheorem{proposition}[theorem]{Proposition}
\newtheorem{corollary}[theorem]{Corollary}

\theoremstyle{definition}
\newtheorem{definition}[theorem]{Definition}

\newtheorem{example}[theorem]{Example}

\newtheorem*{conjecture*}{Conjecture}

\newtheorem{defprop}[theorem]{Definition/Proposition}
\newtheorem{exdef}[theorem]{Example/Definition}

\theoremstyle{remark}
\newtheorem{remark}[theorem]{Remark}

\begin{document}

\title[CSM classes of open quiver loci]{Three formulas
 for CSM classes of open quiver loci}

\author[M. Elkin]{Moriah Elkin}
\address{Cornell University, Ithaca, NY}
\email{mhe45@cornell.edu}

\thanks{}  

\begin{abstract}
In the space of equioriented type $A$ quiver representations, we define subvarieties called ``open quiver loci" by placing strict rank conditions on the maps within representations. The closures of these subvarieties are the quiver loci, whose equivariant cohomology classes are the quiver polynomials of Buch and Fulton. We present one geometric formula and two combinatorial formulas that compute equivariant Chern--Schwartz--MacPherson (CSM) classes of open quiver loci; these classes refine the data of the quiver polynomials. The second combinatorial formula we discuss is in terms of ``chained generic pipe dreams," which modify the pipe dreams of Bergeron and Billey to more strongly resemble the lacing diagrams of Abeasis and Del Fra.
\end{abstract}


\maketitle

\tableofcontents


\section{Introduction}\label{sec:intro}
Historically, much attention has been paid to the ``degeneracy locus" associated to a given map of manifolds: the set of points where the induced map on tangent bundles is below a certain rank. Thom \cite{Thom55} gave a formula for the cohomology class Poincar{\'e} dual to a degeneracy locus, and Porteous \cite{Port} extended the theory to maps of arbitrary vector bundles. Buch and Fulton \cite{BF99} gave a formula in the case of $n$-fold compositions of vector bundle maps over a manifold, where the degeneracy loci are given by forcing each composition of any number of the maps to have a fixed maximal rank. The corresponding cohomology classes are polynomials in the Chern roots of the given vector bundles, and Knutson, Miller, and Shimozono \cite{kms06} gave four combinatorial formulas to compute these polynomials. Later Buch \cite{B04} gave a version of their \emph{pipe dream formula} that removes some cancellation. We present an alternative proof of this formula in Section \ref{sec:polys}, via a new version of their \emph{ratio formula}, that foreshadows proofs in later sections.

The language of equivariant cohomology enables us to work with maps of vector spaces -- type $A$ equioriented quiver representations -- rather than vector bundles. In this setting, degeneracy loci translate to \textit{quiver loci}. Their equivariant cohomology classes, with respect to the group action that changes bases in each of the vector spaces, are called \textit{quiver polynomials}. 

Given a fixed sequence of vector spaces, the quiver loci are closed subvarieties of the corresponding space of quiver representations, and they are given by weak rank conditions: by restricting the maximal allowed rank of each composition of maps. So far, less attention has been paid to the finer objects given by strict rank conditions: by requiring that each composition of maps have some particular rank. We term these subvarieties \textbf{open quiver loci}. They are the orbits of the change-of-basis group action mentioned earlier, and their closures are exactly the quiver loci.

The technology of \textit{equivariant Chern--Schwartz--MacPherson (CSM) classes} (\cite{DG}, \cite{Mac74}, \cite{ohm06}) enables us to associate elements of equivariant cohomology to the \emph{open} quiver loci in a natural way. These classes include the data of the quiver polynomials, as well as the topological Euler characteristic of the open quiver loci. Rim{\'a}nyi \cite{R20} used equivariant localization to compute these classes for general Dynkin type quivers; in this paper, we restrict to type $A$ and are able to give explicit combinatorial formulas.

The first main result of this paper is Theorem \ref{thm:openratio}, which gives a formula for the CSM classes of open quiver loci in terms of (known) CSM classes of Schubert cells in the complete flag variety. We then translate this result into a second formula in terms of the pipe dreams of Bergeron and Billey \cite{BB93}, Proposition \ref{prop:pdopen}.

To improve on the pipe dream formula, we define \textit{chained generic pipe dreams (CGPDs)}, analogues of pipe dreams that closely resemble the \textit{lacing diagrams} of Abeasis and Del Fra \cite{AD85}, and take inspiration from the \textit{generic pipe dreams} of Knutson and Zinn-Justin \cite{kzj25}. Our second main result is Theorem \ref{thm:cgpdopen}, which applies a formula of Su \cite{csu} to express CSM classes of open quiver loci as sums of weights of CGPDs. From Theorem \ref{thm:cgpdopen} we deduce Corollary \ref{cor:cgpd}, a CGPD formula for quiver polynomials.\\

\noindent\textbf{Acknowledgements.} The package I have used for rendering pipe dreams in LaTeX was written by Paul Zinn-Justin. Raj Gandhi, Riccardo Mattarei, and Gabe Udell provided helpful feedback on the exposition, and Anders Buch brought important references to my attention. I am grateful to my advisor Allen Knutson for introducing me to this field and guiding me during this project.

\section{Quiver preliminaries}\label{sec:prelim}
\subsection{Quiver loci and quiver polynomials}\label{subsec:qdefs}
In this section we define and give background on the main objects discussed in the paper, following \cite{kms06}. 

Consider the equioriented type $A_{n+1}$ quiver $\mathrel{\overset{\makebox[0pt]{\mbox{\normalfont\tiny 0}}}{\bullet}} \to \mathrel{\overset{\makebox[0pt]{\mbox{\normalfont\tiny 1}}}{\bullet}} \to  \cdots \to \mathrel{\overset{\makebox[0pt]{\mbox{\tiny $n$}}}{\bullet}}.$
Fix a sequence $V=(V_0,V_1,\dots,V_n)$ of vector spaces over $\C$, and define the variety
\[Hom \coloneqq \Hom(V_0,V_1) \times \cdots \times \Hom(V_{n-1},V_n) =\{V_0 \xrightarrow{\phi_1} V_1 \xrightarrow{\phi_2} \cdots \xrightarrow{\phi_n} V_n\}\] 
of quiver representations on $V$. Given an array of nonnegative integers $\mathbf{r}=(r_{ij})_{0 \leq i\leq j \leq n}$, the \textbf{open quiver locus} $\Omega_{\mathbf{r}}^\circ$ is the set of quiver representations such that $V_i \to V_j$ has rank \textit{exactly} $r_{ij}$.\footnote{We use ``open quiver locus" in the sense of ``open Richardson variety," to indicate that we are not taking the closure of the subvariety of $Hom$ defined by the strict rank conditions. Only the top-dimensonal open quiver locus is truly open in $Hom$, though all are locally closed.} The closure of $\Omega_{\mathbf{r}}^\circ$ is the \textbf{quiver locus} $\Omega_{\mathbf{r}}$, which can be described as the subset of quiver representations whose composite maps $V_i \to V_j$ have rank \textit{at most} $r_{ij}$ for all $i<j$.

Note that there is a natural scheme structure on $\Omega_{\mathbf{r}}$, defined as follows. Fixing a basis for each vector space $V_i$ allows us to express elements of $V_i$ as \textit{row} vectors of length $r_i \coloneqq \dim(V_i)$, and maps $\phi_i$ as matrices with dimensions $r_{i-1} \times r_i$, so the coordinate ring $\C[Hom]$ is a polynomial ring in the ``matrix entry" variables $\{f_{\alpha \beta}^i\}$ where $1\leq i \leq n$, $1 \leq \alpha \leq r_{i-1},$ and $1 \leq \beta \leq r_{i}$. Now the quiver locus is the zero scheme of the ideal
\[I_{\mathbf{r}}=\langle \textrm{minors of size }(1+r_{ij})\textrm{ in } \phi_{i+1}\circ \cdots \circ \phi_j \textrm{ for } i<j \rangle.\]
Lakshmibai and Magyar \cite{lm98} proved that this scheme is reduced using the Zelevinsky map \cite{Zel85}, which we describe in Section \ref{sec:zel}.

To see additional structure of $\Omega_{\mathbf{r}}^\circ$ and $\Omega_{\mathbf{r}}$, consider the action on $Hom$ of the group
\[GL^2=GL(V_0)^2 \times GL(V_1)^2 \times \cdots \times GL(V_{n-1})^2 \times GL(V_n)^2,\]
defined as follows: for each $i$, the first copy of $GL(V_i)$ acts by inverse multiplication on the right of $\Hom(V_{i-1},V_i)$, and the second copy acts by multiplication on the left of $\Hom(V_{i},V_{i+1})$. This action on $\Hom$ can be thought of as ``changing bases of each source and target independently." Since our context involves maps between a fixed set of vector spaces $V_0,\dots,V_n$, it is natural to consider the diagonal embedding
\[GL=GL(V_0)\times \cdots \times GL(V_n) \hookrightarrow GL^2.\]
Now an element $\gamma=(\gamma_0,\dots,\gamma_n) \in GL$ acts on an element $\phi=(\phi_1,\dots,\phi_n)$ via
\[\gamma \cdot \phi=(\gamma_0\phi_1\gamma_1^{-1},\gamma_1\phi_2\gamma_2^{-1},\dots,\gamma_{n-1}\phi_n\gamma_n^{-1}),\]
i.e. by change-of-basis in each vector space $V_i$. This action preserves ranks, so it fixes $\Omega_{\mathbf{r}}^\circ$ and $\Omega_{\mathbf{r}}$.
The goal of Section \ref{sec:polys} is to provide (novel) formulas for \textbf{quiver polynomials}, which we define to be the equivariant cohomology classes $[\Omega_{\mathbf{r}}] \in H^*_{GL}(Hom)$ of the quiver loci.\footnote{Throughout, we will implicitly use Poincar{\'e} duality to identify the equivariant Borel-Moore homology classes defined by $GL$-invariant subvarieties of $Hom$ with the corresponding classes in equivariant cohomology.} Note that $Hom$ is contractible, so $H^*_{GL}(Hom)$ is isomorphic to the $GL$-equivariant cohomology ring of a point, denoted $H^*_{GL}(pt)$.

\subsection{Laces}\label{sec:lace}
Two quiver representations are defined to be \emph{isomorphic} if they are in the same $GL$-orbit. We can define direct sums of quiver representations by taking direct sums of the vector spaces and maps in the natural way, and we call those quiver representations that are not isomorphic to a nontrival direct sum \emph{indecomposable}. The following is a standard result in algebraic quiver theory, c.f. \cite[Sec. 1.1]{lm98}:
\begin{defprop}[\cite{Gabriel}]\label{prop:lace}
Where $0 \leq p \leq q \leq n$, define a quiver representation 
\[I_{p,q}:~\mathrel{\overset{\makebox[0pt]{\mbox{\normalfont\tiny 0}}}{0}} ~\to \cdots \to ~\mathrel{\overset{\makebox[0pt]{\mbox{\tiny p-1}}}{0}} ~\to~ \mathrel{\overset{\makebox[0pt]{\mbox{\tiny p}}}{\raisebox{-1.75pt}{\includegraphics[scale=0.5]{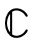}}}} ~\xrightarrow{\sim}\cdots \xrightarrow{\sim} ~\mathrel{\overset{\makebox[0pt]{\mbox{\tiny q}}}{\raisebox{-1.75pt}{\includegraphics[scale=0.5]{Figures/C.png}}}} ~\to ~\mathrel{\overset{\makebox[0pt]{\mbox{\tiny q+1}}}{0}} ~\to \cdots \to ~\mathrel{\overset{\makebox[0pt]{\mbox{\tiny n}}}{0}} \]
to have copies of the field $\C$ at indices $p,\dots,q$ with identity maps between them, and zero vector spaces elsewhere.
Now the $I_{p,q}$ give a complete classification of isomorphism classes of indecomposable quiver representations. In other words, the $I_{p,q}$ are indecomposable, and every $GL$-orbit contains exactly one representation of the form $\bigoplus (I_{p,q})^{s_{pq}}$ for some \textbf{lace array} $\mathbf{s}=(s_{pq})_{0 \leq p \leq q \leq n}$.
\end{defprop}

Given a rank array $\mathbf{r}$ associated to a quiver representation $\phi$, we can recover the corresponding lace array $\mathbf{s}$ by inverting the equation
\begin{equation}
    r_{ij} = \sum_{\substack{p \leq i \\ j \leq q}} s_{pq}, ~~~~~~~ i \leq j
\end{equation}
to arrive at
\begin{equation}
    s_{pq} = r_{pq}-r_{p-1,q}-r_{p,q+1}+r_{p-1,q+1}, ~~~~~~~ p \leq q
\end{equation}
(where we define $r_{ij} \coloneqq 0$ if $i$ and $j$ do not both lie between 0 and $n$; see \cite{kms06}). Therefore $\mathbf{r}$ and $\mathbf{s}$ are equivalent data, and in particular it follows from Proposition \ref{prop:lace} that the open quiver loci $\Omega^\circ_{\mathbf{r}}$ are exactly the $GL$-orbits.

We can picture the $I_{p,q}$ as ``laces" stretching from vertex $p$ of the quiver to vertex $q$. Given a representation $\phi$, the corresponding value of $s_{pq}$ is the number of these ``laces" that appear in a representation isomorphic to $\phi$. Building on this idea, we can represent the full lace array $\mathbf{s}$ pictorially using a graph in the plane called a \textbf{lacing diagram}; these diagrams were developed in \cite{kms06}, derived from those in \cite{AD85}.
The vertices of the graph are arranged in $n+1$ rows with $r_i$ right-justified vertices in row $i$, labeled top-to-bottom. Each vertex in row $i$ may be connected by an edge to at most two other vertices: at most one in row $i-1$, and at most one in row $i+1$. Finally, for every $0 \leq p\leq q\leq n$ there must be exactly $s_{pq}$ distinct connected components with vertices in exactly rows $p,p+1,\dots,q$, corresponding to the summands $I_{p,q}$. A \textbf{lace} is a connected component of a lacing diagram, and an $\mathbf{(i,j)}$-\textbf{lace} is a lace that begins in row $i$ and ends in row $j$.
\begin{example}\label{ex:lace}
Consider the $A_4$ quiver $\bullet \to \bullet \to \bullet \to \bullet$ with sequence of vector spaces $(\C^4, \C^3, \C^3, \C^2)$. The following rank array, lace array, and lacing diagram all describe the same open quiver locus:
\begin{center}
\textbf{r}\,=
\begin{tabular}{c c c c|c}
  3 & 2 & 1 & 0 & $i \diagup j$\\\hline
    &   &   & 4 & 0\\
    &   & 3 & 2 & 1\\
    & 3 & 2 & 1 & 2\\
  2 & 1 & 0 & 0 & 3
\end{tabular}, \hspace{1cm}
\textbf{s}\,=
\begin{tabular}{c c c c|c}
  3 & 2 & 1 & 0 & $i \diagup j$\\\hline
    &   &   & 2 & 0\\
    &   & 0 & 1 & 1\\
    & 0 & 1 & 1 & 2\\
  1 & 1 & 0 & 0 & 3
\end{tabular}, \hspace{1cm}
\begin{tikzpicture}[scale=0.5, baseline={([yshift=-.5ex]current bounding box.center)}]
\node at (3.5,2) {\small 0};
\node at (3.5,-1) {\small 3};
\node at (3.5,0) {\small 2};
\node at (3.5,1) {\small 1};
\filldraw [black] (0,2) circle (4pt);
\filldraw [black] (1,2) circle (4pt);
\filldraw [black] (2,2) circle (4pt);
\filldraw [black] (3,2) circle (4pt);
\filldraw [black] (3,1) circle (4pt);
\filldraw [black] (2,1) circle (4pt);
\filldraw [black] (1,1) circle (4pt);
\filldraw [black] (2,0) circle (4pt);
\filldraw [black] (3,0) circle (4pt);
\filldraw [black] (1,0) circle (4pt);
\filldraw [black] (3,-1) circle (4pt);
\filldraw [black] (2,-1) circle (4pt);
\draw[very thick] (2,1)--(2,2);
\draw[very thick] (3,2)--(3,1)--(3,0);
\draw[very thick] (2,0)--(1,1);
\draw[very thick] (2,-1)--(1,0);
\end{tikzpicture}.
\end{center}
\end{example}
\begin{remark}
There are two notable differences between our lacing diagrams and those discussed in \cite{kms06}. First, we consider lacing diagrams only up to permutation within rows, whereas in \cite{kms06} the order of the basis vectors is fixed at each vertex. Our choice of conventions causes the data of a lacing diagram to be equivalent to that of a lace array, so in particular each open quiver locus has a unique associated lacing diagram. Second, our lacing diagrams are read top-to-bottom rather than left-to-right; the reason for this choice of convention will become clear in Section \ref{sec:opencsm}.
\end{remark}

\subsection{Schubert varieties and the Zelevinsky map}\label{sec:zel}
In this section we describe a useful technique that connects quiver loci to Schubert varieties.

Fix $V=(V_0,V_1,\dots,V_n)$, and write $d=r_0+\dots+r_n$, where $r_i=\dim(V_i)$ as in Section \ref{subsec:qdefs}. Let $B_- \subset GL_d$ be the Borel subgroup of lower triangular matrices. Now $B_-$ acts by left multiplication on $GL_d$, and the quotient $B_- \backslash GL_d$ is the variety of complete flags. Let $B_+$ be the group of upper triangular matrices in $GL_d$, and let $S_d$ be the symmetric group on $d$ letters. Let $v \in S_d$, and identify $v$ with its permutation matrix in $GL_d$. The \textbf{Schubert cell} $X_{v}^\circ \subseteq B_-\backslash GL_d$ is the orbit of $B_- v \in B_-\backslash GL_d$ under the (left) action of $B_+$ by inverse right multiplication. Where $\ell(v)$ is the length of any reduced word for $v$, $X_{v}^\circ$ has codimension $\ell(v)$ in $B_-\backslash GL_d$. Note that we have an equivalent description for $X_{v}^\circ$ in terms of northwest rank conditions: letting $M_{q \times p}$ denote the top $q$ rows and left $p$ columns of a matrix $M$,
\[X_{v}^\circ=B_- \backslash \{M \in GL_d~|~ \rank(M_{q \times p}) = \rank(v_{q \times p}) \textrm{ for all }q,p\}.\]
The \textbf{opposite Schubert cell} $X^{v}_\circ \subset B_-\backslash GL_d$ is the orbit of $B_-v \in B_- \backslash GL_d$ under the (left) action of $B_-$ by inverse right multiplication. It is isomorphic to affine space of dimension $\ell(v)$.
Taking the closures of these orbits gives \textbf{Schubert varieties} $X_v$ and \textbf{opposite Schubert varieties} $X^v$, respectively, and we note that     
\[X_{v}=B_- \backslash \{M \in GL_d~|~ \rank(M_{q \times p}) \leq \rank(v_{q \times p}) \textrm{ for all }q,p\}.\]
For more perspectives on Schubert cells and Schubert varieties, see \cite{bgp25}, \cite{Ful96}, or \cite{gil21}.

We point out a particular opposite Schubert cell that will be important.
Let $w_0$ denote the ``long element" in $S_d$, represented by the antidiagonal permutation matrix. Let $\mathbf{w}_0$ denote the \textit{block} long element with antidiagonal permutation matrices in each diagonal block, where the blocks have sizes $r_0,r_1,\dots,r_n$. Then
\[w_0\mathbf{w}_0=\begin{bmatrix} 0&&&\mathbf{I}_{r_0}\\&&\mathbf{I}_{r_1}&\\&\iddots&&\\\mathbf{I}_{r_n}&&&0\\ \end{bmatrix}, \hspace{1cm} X_\circ^{w_0 \mathbf{w}_0} \cong 
\begin{bmatrix} 
*&*&*&\mathbf{I}_{r_0}\\
*&*&\mathbf{I}_{r_1}&\\
*&\iddots&&\\
\mathbf{I}_{r_n}&&&0\\ \end{bmatrix},\]
where $\mathbf{I}_{r_j}$ denotes the $r_j \times r_j$ identity matrix, and asterisks denote blocks of arbitrary entries.
Note that on the right is a scheme isomorphism identifying $X_\circ^{w_0\mathbf{w}_0}$ with a vector space.

The \textbf{Zelevinsky map} $Z$ was introduced in \cite{Zel85}. We will use the version $Z:Hom \to GL_d$ defined in \cite{kms06} that maps
\[(\phi_1,\phi_2,\dots,\phi_n) \mapsto \begin{bmatrix}
    0 & \cdots & 0 & \phi_1 & \mathbf{I}_{r_0}\\
    0 & \iddots & \phi_2 & \mathbf{I}_{r_1} & 0\\
    0 & \iddots & \mathbf{I}_{r_2} & 0 & 0\\
    \phi_n & \iddots & 0 & 0 & 0\\
    \mathbf{I}_{r_n} & \cdots & 0 & 0 & 0
\end{bmatrix}.\]
The main theorem of \cite{lm98} states that the Zelevinsky map $Z$ induces a scheme isomorphism from each quiver locus $\Omega_{\mathbf{r}}$ to the intersection of a Schubert variety with $X_\circ^{w_0 \mathbf{w}_0}$. In \cite{kms06}, the following definition is made to give the target Schubert variety in terms of $\mathbf{r}$.
\begin{defprop}[Proposition 1.6 and Definition 1.7 of \cite{kms06}]\label{defprop:zel}
Given a rank array $\mathbf{r}$ for $V$ as in Section \ref{subsec:qdefs}, there exists a unique permutation $z(\mathbf{r}) \in S_d$ satisfying the following conditions.
View $z(\mathbf{r})$ as a block matrix, with square blocks of side lengths $r_0,r_1,\dots,r_n$ along the antidiagonal
 from top-right to bottom-left. Now consider the block in the $j$th block row (labeled top to bottom) and $i$th block column (labeled \emph{right to left}):
\begin{enumerate}
    \item If $j \leq i-2$, so the block sits strictly above the main superantidiagonal, then there are no 1 entries in that block.
    \item If $j=i-1$, so the block sits on the main superantidiagonal, then the number of 1 entries in that block equals $r_{j,j+1}$.
    \item If $j \geq i$, so the block sits on or below the main antidiagonal, then the number of 1 entries in that block equals $s_{ij}$.
    \item Within each block row and block column, the 1 entries proceed from northwest to southeast.
\end{enumerate}
This permutation $z(\mathbf{r})$ is the \textbf{Zelevinsky permutation} for the rank array $\mathbf{r}$. 
\end{defprop}
\begin{exdef}\label{exdef:zhom}
For the $A_3$ quiver with $V=(\C^3,\C^3,\C^2)$, let
\[\mathbf{r}_{Hom}=\begin{tabular}{c c c|c}
 2 & 1 & 0 & $i \diagup j$\\\hline
   &   & 3 & 0\\
   & 3 & 3 & 1\\
 2 & 2 & 2 & 2
\end{tabular}, \hspace{1cm}
\textbf{s}_{Hom}=
\begin{tabular}{c c c|c}
2 & 1 & 0 & $i \diagup j$\\\hline
   &   & 0 & 0\\
   & 0 & 1 & 1\\
 0 & 0 & 2 & 2
\end{tabular},\]
so $\Omega_{\mathbf{r}_{Hom}}=Hom$ since all rank conditions are ``full rank."
The associated Zelevinsky permutation matrix is
\begin{center}
   $z(Hom)\coloneqq z(\mathbf{r}_{Hom})=$\hspace{.5cm}
   \begin{NiceTabular}{|cc|ccc|ccc|}[first-row,first-col,last-col]
        & \Block{1-2}{2} & & & 1 & &  & 0& & $i \diagup j$\\\hline
        3 &\cellcolor{gray!30}$\cdot$ &\cellcolor{gray!30} $\cdot$ & 1 & $\cdot$ & $\cdot$ & $\cdot$ & $\cdot$ & $\cdot$ &\\
        4 &\cellcolor{gray!30}$\cdot$ & \cellcolor{gray!30}$\cdot$ & $\cdot$ & 1  & $\cdot$ & $\cdot$ & $\cdot$ & $\cdot$ &0\\
        5 &\cellcolor{gray!30}$\cdot$ & \cellcolor{gray!30}$\cdot$ & $\cdot$ & $\cdot$ & 1 & $\cdot$ & $\cdot$ & $\cdot$ &\\\hline
        1 & 1 & $\cdot$ & $\cdot$ & $\cdot$ & $\cdot$ & $\cdot$ & $\cdot$ & $\cdot$ &\\
        2 &$\cdot$ & 1 & $\cdot$ & $\cdot$ & $\cdot$ & $\cdot$ & $\cdot$ & $\cdot$ &1\\
        6 &$\cdot$ & $\cdot$ & $\cdot$ & $\cdot$ & $\cdot$ & 1 & $\cdot$  & $\cdot$&\\\hline
        7 &$\cdot$ & $\cdot$ & $\cdot$ & $\cdot$ & $\cdot$ & $\cdot$ & 1 & $\cdot$ &\multirow{2}{*}{2}\\
        8 &$\cdot$ & $\cdot$ & $\cdot$ & $\cdot$ & $\cdot$ & $\cdot$ & $\cdot$ & 1 \\\hline
    \end{NiceTabular}
\end{center}
where we have replaced all instances of ``0" with ``$\cdot$" for ease of reading. Here we have also shaded in gray the Rothe diagram, which consists of all matrix entries that are not east or south of any ``1."
We will use the notation $z(Hom)$ in general to mean the Zelevinsky permutation associated to the rank array of the quiver locus $Hom$, in the context of a given quiver and sequence of vector spaces. Note that the Rothe diagram of $z(Hom)$ always consists of exactly the boxes above the block superantidiagonal.
\end{exdef}
\begin{example}
As in Example \ref{ex:lace}, consider the $A_4$ quiver with dimension vector $(4,3,3,2)$, and rank and lace arrays
\[\mathbf{r}=\begin{tabular}{c c c c|c}
  3 & 2 & 1 & 0 & $i \diagup j$\\\hline
    &   &   & 4 & 0\\
    &   & 3 & 2 & 1\\
    & 3 & 2 & 1 & 2\\
  2 & 1 & 0 & 0 & 3
\end{tabular}, \hspace{1cm}
\textbf{s}=
\begin{tabular}{c c c c|c}
  3 & 2 & 1 & 0 & $i \diagup j$\\\hline
    &   &   & 2 & 0\\
    &   & 0 & 1 & 1\\
    & 0 & 1 & 1 & 2\\
  1 & 1 & 0 & 0 & 3
\end{tabular}.\]
The associated Zelevinsky permutation has matrix
\begin{center}
   $z(\mathbf{r})=$\hspace{.5cm}
   \begin{NiceTabular}{|cc|ccc|ccc|cccc|}[first-row,first-col,last-col]
         & \Block{1-2}{3} &  & & 2 & & & 1 & &  \Block{1-4}{0} &&&& $i \diagup j$\\\hline
        6 &\cellcolor{gray!30}$\cdot$ & \cellcolor{gray!30}$\cdot$ & \cellcolor{gray!30}$\cdot$ & \cellcolor{gray!30}$\cdot$ & \cellcolor{gray!30}$\cdot$ & 1 & $\cdot$ & $\cdot$ & $\cdot$ & $\cdot$ & $\cdot$ & $\cdot$ &\multirow{4}{*}{0}\\
        7 &\cellcolor{gray!30}$\cdot$ & \cellcolor{gray!30}$\cdot$ & \cellcolor{gray!30}$\cdot$ & \cellcolor{gray!30}$\cdot$ & \cellcolor{gray!30}$\cdot$ & $\cdot$ & 1 & $\cdot$ & $\cdot$ & $\cdot$ & $\cdot$ & $\cdot$\\
        9 &\cellcolor{gray!30}$\cdot$ & \cellcolor{gray!30}$\cdot$ & \cellcolor{gray!30}$\cdot$ & \cellcolor{gray!30}$\cdot$ & \cellcolor{gray!30}$\cdot$ & $\cdot$ & $\cdot$ & \cellcolor{gray!30}$\cdot$ & 1 & $\cdot$ & $\cdot$ & $\cdot$ \\
        10 &\cellcolor{gray!30}$\cdot$ & \cellcolor{gray!30}$\cdot$ &\cellcolor{gray!30} $\cdot$ &\cellcolor{gray!30} $\cdot$ & \cellcolor{gray!30}$\cdot$ & $\cdot$ & $\cdot$ & \cellcolor{gray!30}$\cdot$ & $\cdot$ & 1 & $\cdot$ & $\cdot$ \\\hline
        3 &\cellcolor{gray!30}$\cdot$ & \cellcolor{gray!30}$\cdot$ & 1 & $\cdot$ & $\cdot$ & $\cdot$ & $\cdot$ & $\cdot$ & $\cdot$ & $\cdot$ & $\cdot$ & $\cdot$\\
        4 &\cellcolor{gray!30}$\cdot$ & \cellcolor{gray!30}$\cdot$ & $\cdot$ & 1 & $\cdot$ & $\cdot$ & $\cdot$ & $\cdot$ & $\cdot$ & $\cdot$ & $\cdot$ & $\cdot$&1\\
        11 &\cellcolor{gray!30}$\cdot$ & \cellcolor{gray!30}$\cdot$ & $\cdot$ & $\cdot$ & \cellcolor{gray!30}$\cdot$ & $\cdot$ & $\cdot$ & \cellcolor{gray!30}$\cdot$ & $\cdot$ & $\cdot$ & 1 & $\cdot$&\\\hline
        1 &1 & $\cdot$ & $\cdot$ & $\cdot$ & $\cdot$ & $\cdot$ & $\cdot$ & $\cdot$ & $\cdot$ & $\cdot$ & $\cdot$ & $\cdot$&\\
        8 &$\cdot$ & \cellcolor{gray!30}$\cdot$ & $\cdot$ & $\cdot$ & \cellcolor{gray!30}$\cdot$ & $\cdot$ & $\cdot$ & 1 & $\cdot$ & $\cdot$ & $\cdot$ & $\cdot$&2\\
        12 &$\cdot$ & \cellcolor{gray!30}$\cdot$ & $\cdot$ & $\cdot$ & \cellcolor{gray!30}$\cdot$ & $\cdot$ & $\cdot$ & $\cdot$ & $\cdot$ & $\cdot$ & $\cdot$ & 1&\\\hline
        2 &$\cdot$ & 1 & $\cdot$ & $\cdot$ & $\cdot$ & $\cdot$ & $\cdot$ & $\cdot$ & $\cdot$ & $\cdot$ & $\cdot$ & $\cdot$&\multirow{2}{*}{3}\\
        5 &$\cdot$ & $\cdot$ & $\cdot$ & $\cdot$ & 1 & $\cdot$ & $\cdot$ & $\cdot$ & $\cdot$ & $\cdot$ & $\cdot$ & $\cdot$\\\hline
    \end{NiceTabular}
\end{center}
where the Rothe diagram is shaded in gray. In two-line notation, the permutation is
\[z(\mathbf{r})=\begin{pmatrix}
    1&2&3&4&5&6&7&8&9&10&11&12\\
    6&7&9&10&3&4&11&1&8&12&2&5
\end{pmatrix}.\]
\end{example}

We can now state the following important proposition, which follows from Theorem 1.14 of \cite{kms06}:
\begin{proposition}\label{prop:zeliso}
The Zelevinsky map $Z$ induces a scheme isomorphism from each quiver locus $\Omega_{\mathbf{r}}$ to the subvariety $X_{z(\mathbf{r})} \cap X_\circ^{w_0\mathbf{w}_0}$ of the flag variety $B_- \backslash GL_d$.
\end{proposition}
The idea of why $z(\mathbf{r})$ is the ``correct" permutation to consider here is as follows.
One should think of $\mathbf{r}$ as defining (weak) \emph{block} northwest rank conditions on the large matrix $Z(\phi_1,\dots,\phi_n)$ for $(\phi_1,\dots,\phi_n) \in \Omega_{\mathbf{r}}$.\footnote{Explicitly, to encode the rank of the map from vertex $i$ to vertex $j$, the rank at the southeast corner of the block in column $i+1$ and row $j-1$ should be at most $r_{ij}+\sum_{k={i+1}}^j r_k$. For instance, the rank of the $i$th superantidiagonal block should be at most $r_{i-1,i}$.
Additionally, all blocks above the block superantidiagonal should have rank 0, and the entire matrix must have full rank.}
A Schubert variety $X_w \subseteq B_- \backslash GL_d$ is instead determined by (weak) northwest rank conditions at \emph{all} matrix entries, with rank jumps at the 1 entries in the matrix for $w$. The Zelevinsky permutation encodes the above block rank conditions into this finer set of rank conditions by jumping rank ``maximally and as soon as possible" within each block. In this way, all matrices with the same or lower northwest rank at the southeast corner of that block are included in $X_{z(\mathbf{r})}$.

\subsection{Double Schubert polynomials and pipe dreams}
We next introduce two important ingredients for the formulas to come in Section \ref{subsec:oldformulas} and beyond.

In \cite{LS82}, Lascoux and Sch{\"u}tzenberger inductively defined the \textit{double Schubert polynomial} $\mathfrak{S}_w(\mathbf{x}-\mathbf{y})$ in the alphabets of variables $\mathbf{x}=(x_1,x_2,\dots,x_d)$ and $\mathbf{y}=(y_1,y_2,\dots,y_d)$ for any permutation $w \in S_d$, $d \in \Z_+$. Double Schubert polynomials represent Schubert classes in equivariant cohomology for the $B_+$-action on the flag variety. They can be computed using diagrams called \emph{pipe dreams}, which were introduced in \cite{BB93} under the name ``rc-graphs."
\begin{definition}
A \textbf{pipe dream} is a $d \times d$ square grid such that each box in the grid is covered by either a \textbf{bump tile} \plaqctr{a} or a \textbf{cross tile} \plaqctr{c}, with exclusively bump tiles on and below the main antidiagonal. Thinking of the lines in the tiled grid as a network of ``pipes," a pipe dream is \textbf{reduced} if each pair of pipes crosses at most once. The set $\mathcal{RP}(v)$ of reduced pipe dreams for a permutation $v \in S_d$ is the set of reduced pipe dreams $D$ such that the pipe entering in row $q$ flows northeast to exit from column $v(q)$.
Given alphabets of weights $\mathbf{x}$ and $\mathbf{y}$, we associate a weight $(\mathbf{x}-\mathbf{y})^D$ to each pipe dream $D$ as follows:
\[(\mathbf{x}-\mathbf{y})^D= \prod_{i,j \in [n]} \begin{cases}
  x_i-y_j & \text{ if } ~\plaqctr{c} ~\text{ at }(i,j)\\[2mm]
  1  & \text{otherwise}
\end{cases}.
\]
\end{definition}
\begin{example}
These are the two pipe dreams for $v=2431$ (here $d=4$).
\begin{center}
\begin{tikzpicture}[scale=0.75]
\node[loop]{\plaq{c}&\plaq{a}&\plaq{a}&\plaq{j}\\\plaq{c}&\plaq{c}&\plaq{j}&\plaq{}\\\plaq{c}&\plaq{j}&\plaq{}&\plaq{}\\\plaq{j}&\plaq{}&\plaq{}&\plaq{}\\};
\node at (-2.35,1.5) {$x_1$};
\node at (-2.35,0.5) {$x_2$};
\node at (-2.35,-0.5) {$x_3$};
\node at (-2.35,-1.5) {$x_4$};
\node at (-1.5,2.35) {$y_1$};
\node at (-0.5,2.35) {$y_2$};
\node at (0.5,2.35) {$y_3$};
\node at (1.5,2.35) {$y_4$};
\end{tikzpicture}
\hspace{1cm}
\begin{tikzpicture}[scale=0.75]\node[loop]{\plaq{c}&\plaq{a}&\plaq{c}&\plaq{j}\\\plaq{c}&\plaq{a}&\plaq{j}&\plaq{}\\\plaq{c}&\plaq{j}&\plaq{}&\plaq{}\\\plaq{j}&\plaq{}&\plaq{}&\plaq{}\\};
\node at (-2.35,1.5) {$x_1$};
\node at (-2.35,0.5) {$x_2$};
\node at (-2.35,-0.5) {$x_3$};
\node at (-2.35,-1.5) {$x_4$};
\node at (-1.5,2.35) {$y_1$};
\node at (-0.5,2.35) {$y_2$};
\node at (0.5,2.35) {$y_3$};
\node at (1.5,2.35) {$y_4$};
\end{tikzpicture}
\end{center}
We have omitted the ``sea" of bumps below the main antidiagonal for ease of reading. For ease of computation, we have labeled the rows and columns of the diagrams with the alphabets of weights $\mathbf{x}$ and $\mathbf{y}$, so that the weight of each cross tile is its row label minus its column label. The weights are 
\[(x_1-y_1)(x_2-y_1)(x_3-y_1)(x_2-y_2)~\text{ and }~(x_1-y_1)(x_2-y_1)(x_3-y_1)(x_1-y_3),\]\
respectively.
\end{example}
\begin{theorem}[\cite{fk96}]\label{thm:pdschubert}
The double Schubert polynomial for $v$ is the sum of the weights of all reduced pipe dreams for $v$:
\[\mathfrak{S}_v(\mathbf{x}-\mathbf{y})=\sum_{D \in \mathcal{RP}(v)}(\mathbf{x}-\mathbf{y})^D.\]
\end{theorem}

\subsection{Knutson-Miller-Shimozono ratio and pipe dream formulas}\label{subsec:oldformulas}
We now introduce two of the formulas given in \cite{kms06} to compute quiver polynomials, to be improved upon in the next section.

Let $T^d=T^{r_0} \times \cdots \times T^{r_n}$ be the diagonal matrices in $GL=GL(V_0) \times \cdots \times GL(V_n)$. In order to compute using pipe dreams, we use the fact that since $Hom$ is contractible, $H^*_{GL}(Hom)$ injects into the polynomial ring $H^*_{T^{d}}(Hom)$, and the torus weights give generators for $H^*_{T^{d}}(Hom) \cong H^*_{T^{d}}(pt)$. We state our pipe dream formula in terms of these weights, with the following indexing conventions. 

In the large $d \times d$ matrix from the Zelevinsky map, let $f^{ij}_{\alpha\beta}$ denote the $(\alpha,\beta)$ entry of the matrix in the $i$th block row (labeled top to bottom) and $j$th block column (labeled \emph{right to left}).
Take $\mathbf{x}=\mathbf{x}^0,\mathbf{x}^1,\dots,\mathbf{x}^n$ to be the (concatenated) alphabet of equivariant parameters for $T^d$ acting on the left. Now via the $T^{2d}$-action coming from the diagonal embedding $GL \hookrightarrow GL^2$ described in Section \ref{subsec:qdefs}, $T^d$ acts on the variable $f_{\alpha \beta}^{i,i+1}$ with weight $x_{\alpha}^{i}-x_\beta^{i+1}$. In the double Schubert polynomial, using $T^d$ in place of $T^{2d}$ corresponds to specializing $\mathbf{y}$ to be a rearrangement of the alphabet $\mathbf{x}$: define $\mathring{\mathbf{x}}\coloneqq \mathbf{x}^n,\mathbf{x}^{n-1},\dots,\mathbf{x}^0$, so the order of the sub-alphabets is reversed, but not the order of the weights within each sub-alphabet. For notational clarity in examples, we rename the sub-alphabets using distinct letters, e.g.
\[\mathbf{x}^0=\mathbf{a}=a_1,a_2,\dots,a_{r_0},\hspace{.5cm}\mathbf{x}^1=\mathbf{b}=b_1,b_2,\dots,b_{r_1},\hspace{.5cm}\mathbf{x}^2=\mathbf{c}=c_1,c_2,\dots,c_{r_2},\hspace{.5cm} \dots.\]

\begin{theorem}[Ratio formula, Theorem 2.9 of \cite{kms06}]\label{thm:oldratio}
The quiver polynomial $[\Omega_{\mathbf{r}}]$ associated to a rank array $\mathbf{r}$ is the ratio
\[[\Omega_{\mathbf{r}}]=\frac{\mathfrak{S}_{z(\mathbf{r})}(\mathbf{x}-\mathring{\mathbf{x}})}{\mathfrak{S}_{z(Hom)}(\mathbf{x}-\mathring{\mathbf{x}})}\]
of specialized double Schubert polynomials.
\end{theorem}
\begin{theorem}[Pipe dream formula, Theorem 5.5 of \cite{kms06}]\label{thm:oldpd}
The quiver polynomial $[\Omega_{\mathbf{r}}]$ associated to a rank array $\mathbf{r}$ is the sum
\[[\Omega_{\mathbf{r}}]=\sum_{D \in \mathcal{RP}(z(\mathbf{r}))} (\mathbf{x}-\mathring{\mathbf{x}})^{D \setminus D_{Hom}},\]
where $(\mathbf{x}-\mathring{\mathbf{x}})^{D \setminus D_{Hom}}$ denotes the weight of the reduced pipe dream $D$ ignoring the weights of the crosses above the block superantidiagonal.
\end{theorem}
\begin{example}[\cite{kms06}]\label{ex:oldpd}
Consider the $A_4$ quiver with $V=(\C^1,\C^3,\C^3,\C^1)$. Below is an example of a rank array $\mathbf{r}$, its corresponding Zelevinsky permutation $z(\mathbf{r})$, and a pipe dream for $z(\mathbf{r})$, labeled with alphabets of weights $\mathbf{x}$ down its rows and $\mathring{\mathbf{x}}$ along its columns. Note the block structure: square blocks of size $r_0=1$, $r_1=3$, $r_2=3$, and $r_3=1$ appear along the antidiagonal, ordered from top-right to bottom-left.
\[
\begin{array}{c@{\hspace{2cm}}c}
\begin{array}{c}
\mathbf{r}=
\begin{tabular}{c c c c|c}
  3 & 2 & 1 & 0 & $i \diagup j$\\\hline
    &   &   & 1 & 0\\
    &   & 3 & 1 & 1\\
    & 3 & 2 & 1 & 2\\
  1 & 1 & 1 & 0 & 3
\end{tabular}\\\\\\
z(\mathbf{r})=\begin{pmatrix}
    1&2&3&4&5&6&7&8\\
    5&2&3&6&1&4&8&7
\end{pmatrix}
\end{array}
&
\begin{tikzpicture}[scale=0.75,baseline={([yshift=-\the\dimexpr\fontdimen22\textfont2\relax]current  bounding  box.center)}]
\node[loop]{
\gplaq{c}&\gplaq{c}&\gplaq{c}&\gplaq{c}&\plaq{a}&\plaq{a}&\plaq{c}&\wplaq{j}\\
\gplaq{c}&\plaq{a}&\plaq{a}&\plaq{a}&\wplaq{a}&\wplaq{a}&\wplaq{j}&\gplaq{}\\
\gplaq{c}&\plaq{a}&\plaq{c}&\plaq{a}&\wplaq{a}&\wplaq{j}&\wplaq{}&\gplaq{}\\
\gplaq{c}&\plaq{a}&\plaq{a}&\plaq{a}&\wplaq{j}&\wplaq{}&\wplaq{}&\gplaq{}\\
\plaq{a}&\wplaq{a}&\wplaq{a}&\wplaq{j}&\gplaq{}&\gplaq{}&\gplaq{}&\gplaq{}\\
\plaq{a}&\wplaq{a}&\wplaq{j}&\wplaq{}&\gplaq{}&\gplaq{}&\gplaq{}&\gplaq{}\\
\plaq{a}&\wplaq{j}&\wplaq{}&\wplaq{}&\gplaq{}&\gplaq{}&\gplaq{}&\gplaq{}\\
\wplaq{j}&\gplaq{}&\gplaq{}&\gplaq{}&\gplaq{}&\gplaq{}&\gplaq{}&\gplaq{}\\
};
\node at (-3.5,4.35) {$d_1$};
\node at (-2.5,4.35) {$c_1$};
\node at (-1.5,4.35) {$c_2$};
\node at (-0.5,4.35) {$c_3$};
\node at (0.5,4.35) {$b_1$};
\node at (1.5,4.35) {$b_2$};
\node at (2.5,4.35) {$b_3$};
\node at (3.5,4.35) {$a_1$};
\node at (-4.35,3.5) {$a_1$};
\node at (-4.35,2.5) {$b_1$};
\node at (-4.35,1.5) {$b_2$};
\node at (-4.35,0.5) {$b_3$};
\node at (-4.35,-0.5) {$c_1$};
\node at (-4.35,-1.5) {$c_2$};
\node at (-4.35,-2.5) {$c_3$};
\node at (-4.35,-3.5) {$d_1$};
\end{tikzpicture}
\end{array}
\]
The pipes crossing at position $(a_1,b_3)$ could cross at any of the seven tiles along the superantidiagonal while still yielding a valid reduced pipe dream for $z(\mathbf{r})$. Independently, moving the cross at position $(b_2,c_2)$ to position $(b_1,c_3)$ or $(b_3,c_1)$ also yields a pipe dream in $\mathcal{RP}(z(\mathbf{r}))$. These operations give all 21 reduced pipe dreams for $z(\mathbf{r})$, so
\[[\Omega_{\mathbf{r}}]=((a_1-b_3)+(b_1-b_2)+(b_2-b_1)+(b_3-c_3)+(c_1-c_2)+(c_2-c_1)+(c_3-d_1))\]\[\cdot ((b_1-c_3)+(b_2-c_2)+(b_3-c_1)).\]
Note that four terms in the first factor cancel, corresponding to twelve ``unnecessary" pipe dreams appearing in the sum. This redundancy will be remedied in the next section.
\end{example}

\section{Streamlined formulas for quiver polynomials}\label{sec:polys}
Here we present ``updated" versions of the ratio and pipe dream formulas in Section \ref{subsec:oldformulas}, which are streamlined in the sense that they have fewer terms; we will use these versions for our results in Section \ref{sec:opencsm}. We begin with the ratio formula. In what follows, the restriction $[X_v\subseteq B_- \backslash GL_d]\Big|_w$ denotes the pullback of the class $[X_v\subseteq B_- \backslash GL_d]$ along the inclusion map of the point $B_-w$ into the flag variety $B_- \backslash GL_d$.
\begin{theorem}[Streamlined ratio formula]\label{thm:newratio}
The quiver polynomial $[\Omega_{\mathbf{r}}]$ associated to a rank array $\mathbf{r}$ is the following ratio of restrictions of equivariant cohomology classes of Schubert varieties:
\[[\Omega_{\mathbf{r}}]=[X_{z(\mathbf{r})} \subseteq B_- \backslash GL_d]\Big|_{w_0\mathbf{w}_0} \Big/ [X_{z(Hom)}\subseteq B_- \backslash GL_d]\Big|_{w_0\mathbf{w}_0}.\]
\end{theorem}
\begin{proof}
Proposition \ref{prop:zeliso} gives the isomorphisms in the following diagram:
\[\begin{tikzcd}[sep=small]
	{\Omega_{\mathbf{r}}} && {X_{z(\mathbf{r})} \cap X_0^{w_0\mathbf{w}_0}} \\
	&&& {X_\circ^{w_0\mathbf{w_0}}} & {B_- \backslash GL_d} \\
	Hom && {X_{z(Hom)} \cap X_\circ^{w_0\mathbf{w}_0}}
	\arrow["\sim", from=1-1, to=1-3]
	\arrow[hook, from=1-1, to=3-1]
	\arrow[hook, from=1-3, to=2-4]
	\arrow[hook, from=1-3, to=3-3]
	\arrow[hook,"\iota", from=2-4, to=2-5]
	\arrow["\sim", from=3-1, to=3-3]
	\arrow[hook, from=3-3, to=2-4]
\end{tikzcd}\]
Now, using the canonical identification $H^*_{T^d}(X) \cong H^*_{T^d}(pt)$ for $X$ contractible, we have the following computation:
\[[\Omega_{\mathbf{r}} \subseteq Hom]=[X_{z(\mathbf{r})} \cap X_\circ^{w_0\mathbf{w}_0} \subseteq X_{z(Hom)} \cap X_\circ^{w_0\mathbf{w}_0}]\]
\[=\frac{[X_{z(\mathbf{r})} \cap X_\circ^{w_0\mathbf{w}_0} \subseteq X_\circ^{w_0\mathbf{w}_0}]}{[X_{z(Hom)} \cap X_\circ^{w_0\mathbf{w}_0} \subseteq X_\circ^{w_0\mathbf{w}_0}]}\]
\[=\frac{\iota^*[X_{z(\mathbf{r})} \subseteq B_- \backslash GL_d]}{\iota^*[X_{z(Hom)}  \subseteq B_- \backslash GL_d]}\]
\[=\frac{[X_{z(\mathbf{r})}\subseteq B_- \backslash GL_d]\Big|_{w_0\mathbf{w}_0}}{[X_{z(Hom)}  \subseteq B_- \backslash GL_d]\Big|_{w_0\mathbf{w}_0}}.\] 
The second equality follows from the axiomatic definition of equivariant cohomology, see for instance \cite[Sec. 1.3]{KZJ07}. The third equality uses the fact that $X_\circ^{w_0\mathbf{w}_0}$ and $X_{z(\mathbf{r})}$ are always transverse.
\end{proof}
We next translate Theorem \ref{thm:newratio} to the following pipe dream formula, which independently follows from the $K$-theoretic statement \cite[Corollary 5.1]{B04} (see also \cite[Corollary 2.2]{BR04}). Theorem \ref{thm:newpd} sums over a smaller set of pipe dreams than Theorem \ref{thm:oldpd}.
\begin{definition}
Given a permutation $v \in S_d$, let $\mathcal{RP}^*(v)$ denote the set of reduced pipe dreams for $v$ that have no crosses on or below the \textbf{block} antidiagonal. (In particular, $\mathcal{RP}^*(v) \subseteq \mathcal{RP}(v)$.)
\end{definition}
\begin{theorem}[Streamlined pipe dream formula]\label{thm:newpd}
The quiver polynomial $[\Omega_{\mathbf{r}}]$ associated to a rank array $\mathbf{r}$ is the sum
   \[[\Omega_{\mathbf{r}}]=\sum_{D \in \mathcal{RP}^*(z(\mathbf{r}))} (\mathbf{x}-\mathring{\mathbf{x}})^{D\setminus D_{Hom}}.\] 
\end{theorem}
\begin{example}
Consider the quiver and rank array from Example \ref{ex:oldpd}. We computed using 21 reduced pipe dreams that the associated quiver polynomial is
\[[\Omega_{\mathbf{r}}]=((a_1-b_3)+(b_1-b_2)+(b_2-b_1)+(b_3-c_3)+(c_1-c_2)+(c_2-c_1)+(c_3-d_1))\]\[\cdot ((b_1-c_3)+(b_2-c_2)+(b_3-c_1)).\]
The pipe dreams associated to the terms that cancel are exactly those with crosses on the block antidiagonal: the new formula sums over only 9 reduced pipe dreams to arrive at the same polynomial,
\[[\Omega_{\mathbf{r}}]=((a_1-b_3)+(b_3-c_3)+(c_3-d_1))\cdot ((b_1-c_3)+(b_2-c_2)+(b_3-c_1)).\]
\end{example}

We will need the following theorem for our proof of Theorem \ref{thm:newpd}.
\begin{theorem}[AJS/Billey formula, \cite{AJS},\cite{Billey}]\label{thm:ajs-b}
Given two permutations $\pi,\rho \in S_n$, let $Q=\sigma_1\sigma_2\dots \sigma_m$ be any reduced word for $\rho$, let $\alpha^i$ be the simple root associated to the simple reflection $\sigma_i$, and let $\beta_i \coloneqq \sigma_1\sigma_2\dots \sigma_{i-1}(\alpha^i)$.\footnote{The superscript on $\alpha^i$ is intended as a reminder that $i$ refers to the index of $\sigma_i$ in $Q$, not the usual index of the simple root -- i.e. $\alpha^i$ will likely not be equal to $\alpha_i$, the simple root corresponding to the reflection $s_i$.} Define
$R(Q,\pi)$ to be the set of subwords of $Q$ that are reduced words for $\pi$. Now in torus-equivariant cohomology of the flag variety $B_- \backslash GL_n$,
\[[X_\pi]\Big|_{\rho}=\sum_{\substack{R=\sigma_{i_1}\sigma_{i_2}\dots \sigma_{i_k}\\ \in R(Q,\pi)}} \beta_{i_1}\beta_{i_2}\cdots\beta_{i_k}.\]
\end{theorem}
\begin{example}
Consider the permutations $\rho=321$ and $\pi=213$ in $S_3$ (written in one-line notation). Using the reduced word $Q=s_1s_2s_1$ for $\rho$, $R(Q,\pi)=\{s_1--,--s_1\}=\{\sigma_1--,--\sigma_3\}$, and Theorem \ref{thm:ajs-b} gives
\[[X_\pi]\Big|_{\rho}=\alpha^1+s_1s_2(\alpha^3)=(y_1-y_2)+(y_2-y_3)=y_1-y_3.\]
If we instead use the reduced word $Q=s_2s_1s_2$, $R(Q,\pi)=\{-s_1-\}=\{-\sigma_2-\}$, and we get the equivalent expression
\[[X_\pi]\Big|_{\rho}=s_2(\alpha^2)=y_1-y_3.\]
\end{example}
\begin{proof}[Proof of Theorem \ref{thm:newpd}] We compute the ratio in Theorem \ref{thm:newratio} using two applications of Theorem \ref{thm:ajs-b}. In both, $\rho=w_0\mathbf{w}_0$. There is exactly one reduced pipe dream for $\rho$, and it has exclusively cross tiles above the block antidiagonal and exclusively bump tiles on and below it, as in the following example for $A_4$ with $V=(\C^1,\C^3,\C^3,\C^1)$ (with $\rho$ written along the left side):
\begin{center}
\begin{tikzpicture}[scale=0.75]
\node[loop]{
\plaq{c}&\plaq{c}&\plaq{c}&\plaq{c}&\plaq{c}&\plaq{c}&\plaq{c}&\wplaq{j}\\
\plaq{c}&\plaq{c}&\plaq{c}&\plaq{c}&\wplaq{a}&\wplaq{a}&\wplaq{j}&\gplaq{}\\
\plaq{c}&\plaq{c}&\plaq{c}&\plaq{c}&\wplaq{a}&\wplaq{j}&\wplaq{}&\gplaq{}\\
\plaq{c}&\plaq{c}&\plaq{c}&\plaq{c}&\wplaq{j}&\wplaq{}&\wplaq{}&\gplaq{}\\
\plaq{c}&\wplaq{a}&\wplaq{a}&\wplaq{j}&\gplaq{}&\gplaq{}&\gplaq{}&\gplaq{}\\
\plaq{c}&\wplaq{a}&\wplaq{j}&\wplaq{}&\gplaq{}&\gplaq{}&\gplaq{}&\gplaq{}\\
\plaq{c}&\wplaq{j}&\wplaq{}&\wplaq{}&\gplaq{}&\gplaq{}&\gplaq{}&\gplaq{}\\
\wplaq{j}&\gplaq{}&\gplaq{}&\gplaq{}&\gplaq{}&\gplaq{}&\gplaq{}&\gplaq{}\\
};
\node at (-4.35,3.5) {$8$};
\node at (-4.35,2.5) {$5$};
\node at (-4.35,1.5) {$6$};
\node at (-4.35,0.5) {$7$};
\node at (-4.35,-0.5) {$2$};
\node at (-4.35,-1.5) {$3$};
\node at (-4.35,-2.5) {$4$};
\node at (-4.35,-3.5) {$1$};
\end{tikzpicture}
\end{center}
The crosses correspond to the simple transpositions in the (unique up to commutation moves) reduced word for $w_0\mathbf{w}_0$. To see this, enumerate all of the cross tiles in the diagram starting in the top-right corner and moving towards to the bottom-left corner; by convention, we choose to move from top-left to bottom-right within each diagonal along the way. Consider the tiles one by one in this order, and record from right to left the positions of the adjacent pipes that were swapped.\footnote{We record right-to-left because we are using the convention where composition of permutations starts from the right, and the correspondence between cross tiles and simple transpositions is immediate when looking top-right to bottom-left.
} Note that these positions are the same along southwest-northeast diagonals. For example, using the $A_2$ quiver and $V=(\C^2,\C^2)$, one can read the word $\textcolor{WildStrawberry}{s_2}\textcolor{Cerulean}{s_3}\textcolor{ForestGreen}{s_1}\textcolor{Brown}{s_2}$ from the following pipe dream for $w_0\mathbf{w}_0$:
\[\begin{tikzpicture}[scale=0.75]
\node[loop]{
\ecplaq{c}{ForestGreen}{lightgray}&\ecplaq{c}{Brown}&\wplaq{a}&\wplaq{j}\\
\ecplaq{c}{WildStrawberry}&\ecplaq{c}{Cerulean}&\wplaq{j}&\wplaq{}\\
\wplaq{a}&\wplaq{j}&\gplaq{}&\gplaq{}\\
\wplaq{j}&\wplaq{}&\gplaq{}&\gplaq{}\\
};
\node at (-2.35,1.5) {$3$};
\node at (-2.35,0.5) {$4$};
\node at (-2.35,-0.5) {$1$};
\node at (-2.35,-1.5) {$2$};
\draw[black,dotted] (-2,1)--(-1,2)--(-1,2.75);
\node at (-1.25,2.35) {\small$s_1$};
\draw[black,dotted] (-2,0)--(0,2)--(0,2.75);
\node at (-0.25,2.35) {\small$s_2$};
\draw[black,dotted] (-2,-1)--(1,2)--(1,2.75);
\node at (0.75,2.35) {\small$s_3$};
\end{tikzpicture}\]

Now picking a subword of a word $Q$ for $\rho$ corresponds to picking certain cross tiles from those present and replacing the rest with bump tiles: i.e., subwords of $\rho$ correspond to pipe dreams with crosses exclusively above the block antidiagonal. It is not hard to see that any reduced pipe dream for a permutation $\pi$ must have exactly $\ell(\pi)$ cross tiles, so we have a correspondence between elements of $R(Q,\pi)$ and elements of $\mathcal{RP}^*(\pi)$. It now only remains to show that the weight of the pipe dream $D$ corresponding to a subword $R=\sigma_{i_1}\sigma_{i_2}\dots \sigma_{i_k}$ equals $\beta_{i_1}\beta_{i_2}\cdots\beta_{i_k}$. To do so, we introduce an illustrative example.

Consider the $A_3$ quiver with $V=(\C^2,\C^2,\C^2)$. Here a word $Q$ for $w_0\mathbf{w}_0$ is $s_4s_5s_3s_4s_2s_3s_1s_4s_2s_5s_3s_4$. Consider the following pipe dream for the subword $R=---s_4--s_1--s_5s_3-=\sigma_4\sigma_7\sigma_{10}\sigma_{11}$ of $Q$, which has weight $(b_1-c_2)(a_1-c_1)(a_2-b_2)(a_1-b_1)$:
\[\begin{tikzpicture}[scale=0.75]
\node[loop]{
\plaq{c}&\plaq{a}&\plaq{c}&\plaq{a}&\wplaq{a}&\wplaq{j}\\
\plaq{a}&\plaq{a}&\plaq{a}&\plaq{c}&\wplaq{j}&\wplaq{}\\
\plaq{a}&\plaq{c}&\wplaq{a}&\wplaq{j}&\gplaq{}&\gplaq{}\\
\plaq{a}&\plaq{a}&\wplaq{j}&\wplaq{}&\gplaq{}&\gplaq{}\\
\wplaq{a}&\wplaq{j}&\gplaq{}&\gplaq{}&\gplaq{}&\gplaq{}\\
\wplaq{j}&\wplaq{}&\gplaq{}&\gplaq{}&\gplaq{}&\gplaq{}\\
};
\node at (-3.35,2.5) {$2$};
\node at (-3.35,1.5) {$1$};
\node at (-3.35,0.5) {$4$};
\node at (-3.35,-0.5) {$6$};
\node at (-3.35,-1.5) {$3$};
\node at (-3.35,-2.5) {$5$};
\draw[black,dotted] (-3,2)--(-2,3)--(-2,3.75);
\node at (-2.25,3.35) {\small$s_1$};
\draw[black,dotted] (-3,1)--(-1,3)--(-1,3.75);
\node at (-1.25,3.35) {\small$s_2$};
\draw[black,dotted] (-3,0)--(0,3)--(0,3.75);
\node at (-0.25,3.35) {\small$s_3$};
\draw[black,dotted] (-3,-1)--(1,3)--(1,3.75);
\node at (0.75,3.35) {\small$s_4$};
\draw[black,dotted] (-3,-2)--(2,3)--(2,3.75);
\node at (1.75,3.35) {\small$s_5$};
\end{tikzpicture}\]
Meanwhile, the summand corresponding to $R$ in the AJS/Billey formula is the product of  $\beta_4=s_4s_5s_3(b_2-c_1)=b_1-c_2$, $\beta_7=s_4s_5s_3s_4s_2s_3(a_1-a_2)=a_1-c_1$, $\beta_{10}=s_4s_5s_3s_4s_2s_3s_1s_4s_2(c_1-c_2)=a_2-b_2$, and $\beta_{11}=s_4s_5s_3s_4s_2s_3s_1s_4s_2s_5(b_1-b_2)=a_1-b_1$. To see why this lines up, consider the marked-up pipe dreams below.

\begin{center}
\includegraphics{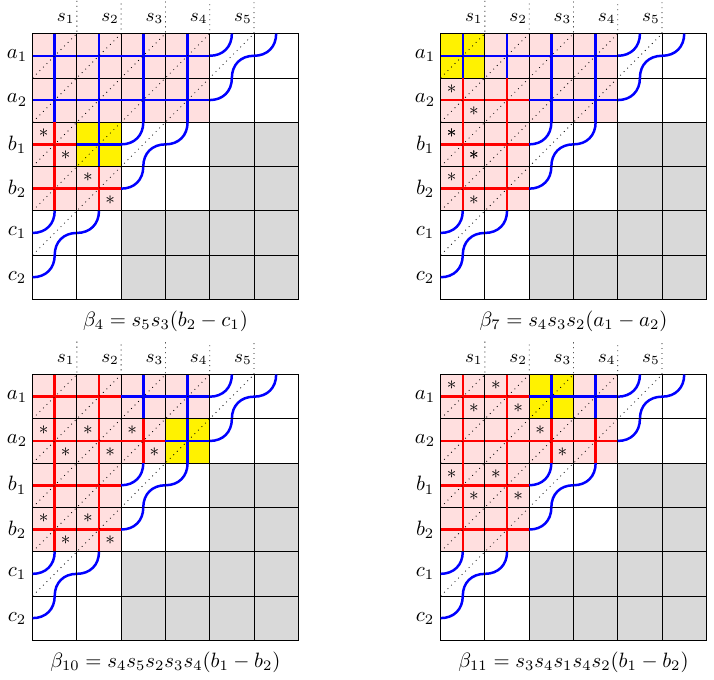}
\end{center}

In each pipe dream, the tiles with red crosses correspond to simple reflections involved in the computation of the respective $\beta_i$, and the simple reflections corresponding to tiles marked with stars are those that have an effect on the calculation; only those starred simple reflections are written in the expressions below each diagram. The starred cross tiles must always appear in a ``hook" left and down from the highlighted tile, where the bottom part of the hook ``reflects" left again along the pipes after reaching the end of the cross tiles. To see this, it is helpful to view the diagram as cut into triangles by the diagonals for the simple reflections. Each triangle corresponds to a weight by following the diagonal strip down and left; the simple root corresponding to a square tile is the difference of the weights of its upper and lower triangle. Now applying the simple reflections in order, if a triangle touches a line for a given reflection, it will be reflected across, changing its weight. By shifting triangles up and down along the diagonal strips (since this doesn't affect the weights they represent), the ``hook" and its effect on the weights becomes clear: the left side of the hook moves the first weight in the simple root $\alpha^i$ left until it is the row label of the highlighted tile, and the lower part of the hook moves the second weight down until its row label equals its column label, then left so the weight records that new row label. The stars drawn into the diagrams illustrate these paths.
In this way, $\beta_i$ becomes the weight contributed by its cross tile to the pipe dream formula: its row label minus its column label.
\end{proof}

\section{CSM class preliminaries}\label{sec:csm}
To study open quiver loci in Section \ref{sec:opencsm}, we will turn to computing their (equivariant) \textbf{Chern--Schwartz--MacPherson classes}. 
Historically, CSM classes of \emph{constructible functions} were defined recursively in homology by MacPherson \cite{Mac74}, proving an existence conjecture of Grothendieck and Deligne. In special cases, they align with earlier work of Schwartz \cite{schwartz1,schwartz2}, hence the name.
The equivariant CSM classes computed in this paper differ in two main ways from the historical definition: they live in equivariant cohomology rather than homology (the extension to the equivariant setting appears first in \cite{ohm06}), and they are homogenized by considering the additional $\C^\times$-action on $T^*Hom$ that fixes all of $Hom$ while scaling its cotangent fibers with weight $-\hbar$.
CSM classes are important for our purposes because they naturally and nontrivially give elements of $H^*_{GL \times \C^\times}(Hom)$ associated to subvarieties that are not necessarily closed. The $\hbar\to\infty$ limit gives the usual classes in $H^*_{GL \times \C^\times}(Hom)$ associated to the subvarieties' closures. In particular, this means that the CSM class of an open quiver locus includes the data of the corresponding quiver polynomial.

The setup of Grothendieck and Deligne proceeds as follows. A \textbf{constructible function} on a complex algebraic variety $Y$ is a finite linear combination of characteristic functions of closed subvarieties; we let $\mathcal{F}(Y)$ denote the group (under pointwise addition) of such functions. Now a proper morphism $f:Y \to X$ of complex algebraic varieties gives rise to a pushforward map $f_*:\mathcal{F}(Y) \to \mathcal{F}(X)$, by linearly extending the following definition for the characteristic function of a locally closed subvariety $W \subseteq Y$: where $\chi$ is the topological Euler characteristic and $p \in X$, we define $f_*(\mathbb{1}_W)(p)\coloneqq \chi(f^{-1}(p) \cap W)$. With these definitions, $\mathcal{F}$ is in fact a covariant functor from the category of algebraic varieties to the category of abelian groups, analogous to the homology functor $H_*$. MacPherson constructed a natural transformation $c: \mathcal{F} \to H_*$ satisfying the ``normalization" requirement that if $Y$ is non-singular, then $c(\mathbb{1}_Y)$ is Poincar{\'e} dual to the total Chern class of the tangent bundle of $Y$.
Note that the definition in terms of constructible functions necessitates an important and useful property of CSM classes: that for disjoint subvarieties $A,B \subseteq Y$, $\csm(A \subseteq Y)+\csm(B \subseteq Y)=\csm(A \sqcup B \subseteq Y)$.

In the equivariant setting, the variety $Y$ comes with an action of a reductive linear group $G$, any map $f$ must be $G$-equivariant, and we consider only the $G$-invariant constructible functions $\F^G(Y)$ (which are constructed from characteristic functions of locally closed $G$-invariant subvarieties, rather than all locally closed subvarieties). To homogenize in the case that $Y$ is smooth, we additionally consider the $\C^\times$ action that fixes 
$Y$ and scales its cotangent fibers with weight $-\hbar$ (or equivalently, scales its tangent fibers with weight $\hbar$).
Due to \cite{ohm06}, in this setting we have a natural transformation $c_G: \mathcal{F}^{G} \to H^*_{G\times \C^\times}$ such that $c_G(\mathbb{1}_Y)$ is $e_{G\times \C^{\times}}(TY)$, the $(G\times \C^\times)$-equivariant Euler class of the tangent bundle of $Y$ (which is the $\hbar$-homogenization of the $G$-equivariant total Chern class of $TY$).

We will soon set $Y \coloneqq Hom$ and $G \coloneqq GL$, and consider the CSM classes associated to characteristic functions of open quiver loci, though we will slightly abuse notation by writing $\csm(\Omega_{\mathbf{r}}^\circ \subseteq Hom)$ for $c_{GL}(\mathbb{1}_{\Omega_{\mathbf{r}}^\circ} \subseteq Hom)$. 
First, for completeness, we present two widely known lemmas. The proof we present of Lemma \ref{lem:vs} was communicated to us by Leonardo Mihalcea.
\begin{lemma}\label{lem:vs}
Let $A$ be a smooth closed subvariety of a smooth variety $M$, and let $C_AM$ denote the conormal bundle to $A$ in $M$. Then $\csm(A \subseteq M)=(-1)^{\dim A} [C_AM \subseteq T^*M]$.
\end{lemma}
\begin{proof}
From the short exact sequence of vector bundles
\begin{center}
    \includegraphics{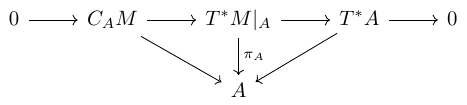}
\end{center}
we conclude that $[C_AM \subseteq T^*M|_A] = \pi^*_A (e(T^*A))$ in equivariant cohomology (see \cite[Example 6.3.5]{fultonIT}).
Additionally, letting $i_A$ and $i_M$ be zero sections, the following diagram is a pullback square:
\begin{center}
    \includegraphics{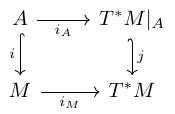}
\end{center}
Therefore for any $[X] \in H^*_{G\times \C^\times}(T^*M|_A)$, $i_M^* j_* ([X]) = i_* i_A^* ([X])$ (see \cite[Section 3.6]{andfult}).
We now compute that
\begin{align*}
    [C_AM \subseteq T^*M] &= i_M^*j_*[C_AM \subseteq T^*M|_A]\\
    &= i_*i_A^*(\pi^*_A (e (T^*A)))\\
    &= i_*e (T^*A)\\
    &= (-1)^{\dim(A)} (i_* e(TA))\\
    &= (-1)^{\dim(A)} \csm(A \subseteq M).\qedhere
\end{align*}
\end{proof}
\begin{lemma}[c.f. Theorem 4.2 of \cite{amss24}]\label{lem:csm}
Given $A$ and $B$ subvarieties of smooth varieties $M$ and $N$ respectively, we have
\[\csm(A\times B \subseteq M \times N)= \csm(A \subseteq M)\cdot\csm(B \subseteq N)\]
under the K{\"u}nneth map
\[H^*_{G \times \C^\times}(T^*M) \otimes H^*_{G \times \C^\times}(T^*N) \to H^*_{G \times (\C^\times)^2}(T^*M \times T^*N)\]
(where the two parameters for the $(\C^\times)^2$-action are set equal).
\end{lemma}

We conclude this section by computing the CSM classes of two (very small) open quiver loci.
\begin{example}
Consider the $A_2$ quiver $\bullet \to \bullet$ with sequence of vector spaces $(\C^1,\C^1)$. $GL=(GL_1)^2$ acts on $Hom= \{\C \xrightarrow{q} \C\}$ with two orbits, $\{q=0\} \sqcup \{q \neq 0\}$, corresponding to the rank arrays
\[\mathbf{r}_1=
\begin{tabular}{c c|c}
    1 & 0 & $i \diagup j$\\\hline
     & 1 & 0\\
    1 & 0 & 1
\end{tabular}, \hspace{1cm} \mathbf{r}_2=
\begin{tabular}{c c|c}
    1 & 0 & $i \diagup j$\\\hline
     & 1 & 0\\
    1 & 1 & 1
\end{tabular}.\]
We will use the coordinate $p$ for the cotangent vectors, and label the torus weights $x_1-x_2$ on the $q$-axis and $-\hbar-(x_1-x_2)$ on the $p$-axis.

From Lemma \ref{lem:vs}, $\csm(\Omega_{\mathbf{r}_1}^\circ \subseteq Hom)$ is the class of the $p$-axis in $T^*Hom$, with no added sign since $\Omega_{\mathbf{r}_1}^\circ$ has dimension 0. Therefore $\csm(\Omega_{\mathbf{r}_1}^\circ \subseteq Hom)=x_1-x_2$.
To compute $\csm(\Omega_{\mathbf{r}_2}^\circ \subseteq Hom)$ we use additivity, noting that $\Omega_{\mathbf{r}_1}^\circ \sqcup \Omega_{\mathbf{r}_2}^\circ = Hom$. Again from Lemma \ref{lem:vs}, $\csm(Hom \subseteq Hom)$ is the class of the $q$-axis in $T^*Hom$, with an added factor of $-1$ since $Hom$ is one-dimensional. It follows that
\[\csm(\Omega_{\mathbf{r}_2}^\circ \subseteq Hom)=\csm(Hom)-\csm(\Omega_{\mathbf{r}_1}^\circ)=-(-\hbar-x_1+x_2) -(x_1-x_2)=\hbar.\]
\end{example}

\section{Computing CSM classes of open quiver loci}\label{sec:opencsm}
In this section we give three formulas for CSM classes of open quiver loci: one ratio formula, one using pipe dreams, and one using ``chained generic pipe dreams," which we will define. The statements and proofs somewhat resemble those in Section \ref{sec:polys}.

\subsection{Ratio formula}
We begin by establishing an analogue of Proposition \ref{prop:zeliso} for open quiver loci. We can no longer use $z(\mathbf{r})$ to encode a given rank array $\mathbf{r}$ into rank conditions on the matrices in the image of the Zelevinsky map. Instead, we must consider every permutation satisfying the block rank conditions given by $\mathbf{r}$, as described in the discussion after Proposition \ref{prop:zeliso}. 
\begin{definition}
Given a rank array $\mathbf{r}$ for $V$ as in Section \ref{subsec:qdefs}, write $d=r_0+\dots+r_n$ as in Section \ref{sec:zel}. We denote by $\perm(\mathbf{r})$ the set of all
permutations $w \in S_d$
such that, viewing $w$ as a block matrix like $z(\mathbf{r})$, each block of $w$ contains the same number of 1 entries as the corresponding block of $z(\mathbf{r})$.
\end{definition}
\begin{example}
For the quiver $A_2$ with $V=(\C^1,\C^2)$, let $\mathbf{r}=
\begin{tabular}{c c|c}
    1 & 0 & $i \diagup j$\\\hline
     & 1 & 0\\
    2 & 1 & 1
\end{tabular}.$
The associated block rank conditions are
$~\begin{tabular}{c c|c}
    1 & 0 & $i \diagup j$\\\hline
    1 & 1 & 0\\
    2 & 3 & 1
\end{tabular}.$
Now 
\[\perm(\mathbf{r})=\left\{
\left[\begin{tabular}{cc|c}
    1 & 0 & 0 \\\hline
    0 & 1 & 0 \\
    0 & 0 & 1 \\
\end{tabular}\right],
\left[\begin{tabular}{cc|c}
    1 & 0 & 0 \\\hline
    0 & 0 & 1 \\
    0 & 1 & 0 \\
\end{tabular}\right],
\left[\begin{tabular}{cc|c}
    0 & 1 & 0 \\\hline
    1 & 0 & 0 \\
    0 & 0 & 1 \\
\end{tabular}\right],
\left[\begin{tabular}{cc|c}
    0 & 1 & 0 \\\hline
    0 & 0 & 1 \\
    1 & 0 & 0 \\
\end{tabular}\right]
\right\}.\]
\end{example}
\begin{proposition}\label{prop:openzeliso}
The Zelevinsky map $Z$ induces a scheme isomorphism from each open quiver locus $\Omega^\circ_{\mathbf{r}}$ to the disjoint union
\[\bigsqcup_{v \in \perm(\mathbf{r})} \left( X_{v}^\circ \cap X_\circ^{w_0\mathbf{w}_0} \right) \quad \subset B_- \backslash GL_d.\]
\end{proposition}
\begin{proof}
$Z$ gives a scheme isomorphism from $\Omega^\circ_{\mathbf{r}}$ to its image $Z(\Omega^\circ_{\mathbf{r}})$ in $GL_d$, and therefore to the image $Y_{\mathbf{r}} \coloneqq B_-\backslash Z(\Omega^\circ_{\mathbf{r}})$ in $B_-\backslash GL_d$, since each element of $Z(\Omega^\circ_{\mathbf{r}})$ represents a unique $B_-$-orbit. To set-theoretically identify $Y_{\mathbf{r}}$ with $\bigsqcup_{v \in \perm(\mathbf{r})} X_{v}^\circ \cap X_\circ^{w_0\mathbf{w}_0}$, first note that the image of $Z$ is certainly contained in $X_\circ^{w_0\mathbf{w}_0}$. Now as described after Proposition \ref{prop:zeliso}, $\mathbf{r}$ can be viewed as a set of strict block northwest rank conditions on $Z(\Omega^\circ_{\mathbf{r}})$, and these block rank conditions are fulfilled precisely by the block rank jumps appearing in $z(\mathbf{r})$. Translating these block rank conditions into rank conditions for all entries, each element of $Y_{\mathbf{r}}$ is in the Schubert cell for some permutation with the same block rank jumps as $z(\mathbf{r})$, i.e. for some permutation in $\perm(\mathbf{r})$. Since the block rank conditions entirely define $Y_{\mathbf{r}}$ as a subset of $X_\circ^{w_0\mathbf{w}_0}$, we are done.
\end{proof}

We now have the tools to prove the ratio formula.
\begin{theorem}[Ratio formula for open quiver loci]\label{thm:openratio}
Given a rank array $\mathbf{r}$, we have the following expression for the CSM class of the corresponding open quiver locus:
\[\csm(\Omega_{\mathbf{r}}^\circ \subset Hom) =
\left(\sum_{v \in \perm(\mathbf{r})}\csm(X_{v}^\circ \subset B_- \backslash GL_d)|_{w_0\mathbf{w_0}}\right) \Bigg/ [X_{z(Hom)}\subset B_- \backslash GL_d]|_{w_0\mathbf{w_0}}.\]
\end{theorem}
\begin{proof}
Using Proposition \ref{prop:openzeliso} and the additivity of CSM classes,
\[\csm(\Omega_{\mathbf{r}}^\circ \subset Hom)= \sum_{v \in \perm(\mathbf{r})} \csm(X_{v}^\circ \cap X_\circ^{w_0\mathbf{w}_0} \subset X_{z(Hom)} \cap X_\circ^{w_0\mathbf{w}_0}).\]
Recall that
\[
X_\circ^{w_0\mathbf{w}_0} \cong 
\begin{bmatrix} 
*&*&*&\mathbf{I}_{r_0}\\
*&*&\mathbf{I}_{r_1}&\\
*&\iddots&&\\
\mathbf{I}_{r_n}&&&0\\ \end{bmatrix}
\hspace{.75cm}\textrm{ and }\hspace{.75cm}
X_{z(Hom)} \cap X_\circ^{w_0\mathbf{w}_0} \cong \begin{bmatrix} 
0&0&*&\mathbf{I}_{r_0}\\
0&*&\mathbf{I}_{r_1}&\\
*&\iddots&&\\]
\mathbf{I}_{r_n}&&&0\\ \end{bmatrix}
,\]
so $X_\circ^{w_0\mathbf{w}_0} \cong (X_{z(Hom)}\cap X_\circ^{w_0\mathbf{w}_0}) \times W$, where $W$ is the vector space of matrices with entries strictly above the block superantidiagonal. Now from Lemma \ref{lem:csm}, for each $v \in \perm(\mathbf{r})$,
\[\csm(X_{v}^\circ \cap X_\circ^{w_0\mathbf{w}_0} \subset X_{z(Hom)} \cap X_\circ^{w_0\mathbf{w}_0}) = \frac{\csm(X_{v}^\circ \cap X_\circ^{w_0\mathbf{w}_0} \subset X_\circ^{w_0\mathbf{w}_0})}{\csm(0 \subset W)}=\frac{\csm(X_{v}^\circ \subset B_- \backslash GL_d)|_{w_0\mathbf{w_0}}}{\csm(0 \subset W)}.\]
Referencing Lemma \ref{lem:vs}, we compute that (where $\mathrm{wt}(W)$ is the product of the torus weights on $W$)
\[\csm(0 \subset W)=(-1)^{0}[T^*0 \subseteq T^*W] =\text{wt}(W)\]
\[=[X_{z(Hom)} \cap X_\circ^{w_0\mathbf{w}_0}\subset X_\circ^{w_0\mathbf{w}_0}]=[X_{z(Hom)}\subset B_- \backslash GL_d]|_{w_0\mathbf{w_0}}\]
and the claim follows.
\end{proof}

\subsection{Pipe dream formula}\label{sec:csm-pd}
The following theorem translates the ratio formula into the language of pipe dreams.
\begin{proposition}[Pipe dream formula for open quiver loci]\label{prop:pdopen}
For any $v \in S_d$, let $\mathcal{P}(v)$ be the set of (not necessarily reduced!) pipe dreams for $v$, and let $c(D)$ be the set of cross tiles in a given pipe dream $D$, so $|c(D)|$ is the number of cross tiles in $D$.
Given a rank array $\mathbf{r}$, we have the following expression for the CSM class of the corresponding open quiver locus:
\[\csm(\Omega_{\mathbf{r}}^\circ \subset P_- \backslash GL_d) =
\sum_{v \in \perm(\mathbf{r})} 
\sum_{D \in \mathcal{P}(v)} \hbar^{\ell(w_0\mathbf{w}_0)-|c(D)|}(\mathbf{x}-\mathring{\mathbf{x}})^{D\setminus D_{Hom}}.\]
\end{proposition}

The proof is similar to that of Theorem \ref{thm:newpd}. In \cite{csu}, Su gives an analogue of the AJS/Billey formula (Theorem \ref{thm:ajs-b}) for stable envelopes of cotangent bundles of flag varieties. These are linked to characteristic cycles by Maulik and Okounkov \cite{MO}, whose cohomology classes are identified with CSM classes by Ginsburg \cite{g86}, enabling us to apply Su's formula to our calculations.
Many of these pieces are explicitly put together in \cite[Section 9]{AMSS23}.
\begin{theorem}[Corollary 9.8 of \cite{AMSS23}]\label{thm:su}
Given two elements $\pi,\rho \in S_n$, let $Q=\sigma_1\sigma_2\dots \sigma_m$ be any reduced word for $\rho$. Let $\alpha^i$ be the simple root associated to the simple reflection $\sigma_i$, and write $\beta_i \coloneqq \sigma_1\sigma_2\dots \sigma_{i-1}(\alpha^i)$ with $R(Q)\coloneqq\{\beta_i|1 \leq i \leq m\}$. Write $R^+$ for the set of positive roots of $GL_d$.
Define $W(Q,v)$ to be the set of subwords of $Q$ that are (not necessarily reduced!) words for a permutation $v$. Now
\[\csm(X_\pi^\circ \subseteq B_- \backslash GL_d)\Big|_{\rho}=\left(\prod_{\alpha \in R^+ \setminus R(Q)} (\alpha-\hbar)\right)
\sum_{\substack{W=\sigma_{i_1}\sigma_{i_2}\dots \sigma_{i_k}\\ \in W(Q,v)}} \hbar^{m-k}\beta_{i_1}\beta_{i_2}\cdots\beta_{i_k}.\]
\end{theorem}
\begin{remark}
An analogue of Theorem \ref{thm:su} that computes motivic Chern classes is not currently available in the literature. However, such a formula would likely lead to $K$-theoretic analogues of the diagrammatic formulas given in Proposition \ref{prop:pdopen} and Theorem \ref{thm:cgpdopen}.
\end{remark}
\begin{proof}[Proof of Proposition \ref{prop:pdopen}]
First, note that plugging in $\rho=w_0\mathbf{w}_0$ to Theorem \ref{thm:su} causes $R^+ \setminus R(Q)$ to be empty, so the product on the left is 1. Now as in the proof of Theorem \ref{thm:newpd}, each subword $\sigma_{i_1}\sigma_{i_2}\dots \sigma_{i_k}$ of $Q$ corresponds to a pipe dream with cross tiles strictly above the block superantidiagonal, with weight $\beta_{i_1}\beta_{i_2}\cdots\beta_{i_k}$. 
There are two main differences between this formula and Theorem \ref{thm:ajs-b} that require us to adjust the pipe dream formula: first, the subword and therefore the pipe dream is not required to be reduced; and second, we have a factor of $\hbar$ for each bump tile present above the block superantidiagonal. 
Applying this interpretation of Theorem \ref{thm:su} to Theorem \ref{thm:openratio}  (along with the analogous interpretation of Theorem \ref{thm:ajs-b} in the denominator) yields the statement of the proposition.
\end{proof}

\begin{example}\label{ex:a3}
Consider the $A_3$ quiver with $V=(\C^1,\C^2,\C^1)$, and take
$\mathbf{r}=
\begin{tabular}{ccc|c}
    2 & 1 & 0 & $i \diagup j$\\\hline
  &   & 1 &0\\
  &  2 & 1 & 1\\
   1 & 1 & 1 & 2 
\end{tabular}.$
The corresponding block rank conditions are
$\begin{tabular}{ccc|c}
    2 & 1 & 0 & $i \diagup j$\\\hline
 0 &  1 & 1 &0\\
 1 &  3 & 3 & 1\\
   1 & 3& 4 & 2 
\end{tabular}.$
The five pipe dreams for permutations $v \in \perm(\mathbf{r})$ are
\begin{center}
\adjustbox{max width=\textwidth}{
\begin{tikzpicture}[scale=0.75]
\node[loop]{
\gplaq{c}&\plaq{a}&\plaq{a}&\wplaq{j}\\
\plaq{a}&\wplaq{a}&\wplaq{j}&\gplaq{}\\
\plaq{a}&\wplaq{j}&\wplaq{}&\gplaq{}\\
\wplaq{j}&\gplaq{}&\gplaq{}&\gplaq{}\\};
\node at (-2.35,1.5) {$a$};
\node at (-2.35,0.5) {$b_1$};
\node at (-2.35,-0.5) {$b_2$};
\node at (-2.35,-1.5) {$c$};
\node at (-1.5,2.35) {$c$};
\node at (-0.5,2.35) {$b_1$};
\node at (0.5,2.35) {$b_2$};
\node at (1.5,2.35) {$a$};
\node at (0,-2.35) {$v=2134$};
\end{tikzpicture}
\hspace{0.2cm}
\begin{tikzpicture}[scale=0.75]\node[loop]{\gplaq{c}&\plaq{a}&\plaq{c}&\wplaq{j}\\
\plaq{a}&\wplaq{a}&\wplaq{j}&\gplaq{}\\
\plaq{c}&\wplaq{j}&\wplaq{}&\gplaq{}\\
\wplaq{j}&\gplaq{}&\gplaq{}&\gplaq{}\\};
\node at (-1.5,2.35) {$c$};
\node at (-0.5,2.35) {$b_1$};
\node at (0.5,2.35) {$b_2$};
\node at (1.5,2.35) {$a$};
\node at (0,-2.35) {$v=2134$};
\end{tikzpicture}
\hspace{0.2cm}
\begin{tikzpicture}[scale=0.75]
\node[loop]{
\gplaq{c}&\plaq{c}&\plaq{a}&\wplaq{j}\\
\plaq{c}&\wplaq{a}&\wplaq{j}&\gplaq{}\\
\plaq{a}&\wplaq{j}&\wplaq{}&\gplaq{}\\
\wplaq{j}&\gplaq{}&\gplaq{}&\gplaq{}\\};
\node at (-1.5,2.35) {$c$};
\node at (-0.5,2.35) {$b_1$};
\node at (0.5,2.35) {$b_2$};
\node at (1.5,2.35) {$a$};
\node at (0,-2.35) {$v=3214$};
\end{tikzpicture}
\hspace{0.2cm}
\begin{tikzpicture}[scale=0.75]
\node[loop]{
\gplaq{c}&\plaq{c}&\plaq{a}&\wplaq{j}\\
\plaq{a}&\wplaq{a}&\wplaq{j}&\gplaq{}\\
\plaq{a}&\wplaq{j}&\wplaq{}&\gplaq{}\\
\wplaq{j}&\gplaq{}&\gplaq{}&\gplaq{}\\};
\node at (-1.5,2.35) {$c$};
\node at (-0.5,2.35) {$b_1$};
\node at (0.5,2.35) {$b_2$};
\node at (1.5,2.35) {$a$};
\node at (0,-2.35) {$v=3124$};
\end{tikzpicture}
\hspace{0.2cm}
\begin{tikzpicture}[scale=0.75]
\node[loop]{
\gplaq{c}&\plaq{a}&\plaq{a}&\wplaq{j}\\
\plaq{c}&\wplaq{a}&\wplaq{j}&\gplaq{}\\
\plaq{a}&\wplaq{j}&\wplaq{}&\gplaq{}\\
\wplaq{j}&\gplaq{}&\gplaq{}&\gplaq{}\\};
\node at (-1.5,2.35) {$c$};
\node at (-0.5,2.35) {$b_1$};
\node at (0.5,2.35) {$b_2$};
\node at (1.5,2.35) {$a$};
\node at (0,-2.35) {$v=2314$};
\end{tikzpicture}}
\end{center}
and Proposition \ref{prop:pdopen} gives
\[\csm(\Omega_{\mathbf{r}}^\circ \subset {B_-} \backslash GL_d)=\hbar^4+\hbar^2(a-b_2)(b_2-c)+\hbar^2(b_1-c)(a-b_1)+\hbar^3(a-b_1)+\hbar^3(b_1-c).\]
\end{example}

\subsection{Chained generic pipe dream formula}\label{sec:cgpd}
We now define a more natural object for this computation, inspired by the ``generic pipe dreams" of \cite{kzj25}.
\begin{definition}
Fix a sequence of vector spaces $V=(V_0,V_1,\dots,V_n)$ with $\dim(V_i)=r_i$, and a rank array $\mathbf{r}$ with lace array $\mathbf{s}$ for $V$. 
A \textbf{chained generic pipe dream tile} is any of the following, where blue and red represent arbitrary distinct colors:
\begin{center}
\plaqctr{a}\hspace{.6cm}
\begin{tikzpicture}[baseline={([yshift=-\the\dimexpr\fontdimen22\textfont2\relax]current  bounding  box.center)}]
\node[bgplaq,rectangle,draw,minimum size=\loopcellsize,transform shape] (plaq) {};\useasboundingbox;\begin{scope}[x=\loopcellsize,y=\loopcellsize] \draw[/linkpattern/edge,\plaqwest,\plaqnorth,color=red] (0,0.5) .. controls (0,0.2) and (-0.2,0) .. (-0.5,0); \draw[/linkpattern/edge,\plaqeast,\plaqsouth] (0,-0.5) .. controls (0,-0.2) and (0.2,0) .. (0.5,0); \end{scope}\end{tikzpicture}
\hspace{.6cm}
\plaqctr{j}\hspace{.6cm}
\plaqctr{r}\hspace{.6cm}
\begin{tikzpicture}[baseline={([yshift=-\the\dimexpr\fontdimen22\textfont2\relax]current  bounding  box.center)}]
\node[bgplaq,rectangle,draw,minimum size=\loopcellsize,transform shape] (plaq) {};\useasboundingbox;\begin{scope}[x=\loopcellsize,y=\loopcellsize] \draw[/linkpattern/edge,\plaqsouth,\plaqnorth] (0,0.5) -- (0,-0.5);
\draw[/linkpattern/edge,\plaqwest,\plaqeast,color=red] (0.5,0) -- (-0.5,0); \end{scope}
\end{tikzpicture}\hspace{.6cm}
\plaqctr{h}\hspace{.6cm}
\plaqctr{v}\hspace{.6cm}
\plaqctr{}
\end{center}
A \textbf{chained generic pipe dream (CGPD) for $V$} is an arrangement of these tiles into $n+1$ rectangles, joined at their northeast and southwest corners and labeled $0,\dots,n$ from northeast to southwest. Rectangle $i$ has dimensions $r_i \times r_{i+1}$ (where we define the ($n+1$)st rectangle to have dimensions $r_n \times 0$), and $r_i$ pipes enter it along its east side.
If any pipes exit from the $j$th column at the bottom of a rectangle, they continue directly southwest to enter from the $j$th row in the subsequent rectangle. Pipes exiting through the west sides of rectangles permanently leave the diagram, and each pipe is assigned a color corresponding to the last rectangle in which it appears. Two pipes of the same color are not allowed to cross.

A \textbf{CGPD for $\mathbf{r}$} has the additional property
that for any $0 \leq i \leq j \leq n$, the number of pipes appearing in exactly rectangles $i,i+1,\dots,j$ is $s_{ij}$. We denote the set of CGPDs for $\mathbf{r}$ by $\mathcal{CGPD}(\mathbf{r})$.
\end{definition}
\begin{example}\label{ex:cgpd}
Let $V=(\C^1,\C^1,\C^2,\C^3,\C^2,\C^1)$. Consider the following rank array $\mathbf{r}$ and corresponding lacing diagram $L$.
\begin{center}
$\mathbf{r}=
\begin{tabular}{cccccc|c}
5 & 4 & 3 & 2 & 1 & 0 & $i \diagup j$\\\hline
 &  &  &  &  & 1 & 0\\
 &  &  &  & 1 & 1 & 1\\
 &  &  & 2 & 1 & 1 & 2\\
 &  & 3 & 1 & 0 & 0 & 3\\
 & 2 & 1 & 1 & 0 & 0 & 4\\
1 & 1 & 1 & 1 & 0 & 0 & 5\\
\end{tabular}$\hspace{2cm}
$L=$\hspace{0.25cm}\begin{tikzpicture}[scale=0.5, baseline={([yshift=-.5ex]current bounding box.center)}]
\node at (2.5,6) {\small 0};
\node at (2.5,5) {\small 1};
\node at (2.5,4) {\small 2};
\node at (2.5,3) {\small 3};
\node at (2.5,2) {\small 4};
\node at (2.5,1) {\small 5};
\filldraw [black] (2,6) circle (4pt);
\filldraw [black] (2,5) circle (4pt);
\filldraw [black] (2,4) circle (4pt);
\filldraw [black] (1,4) circle (4pt);
\filldraw [black] (0,3) circle (4pt);
\filldraw [black] (1,3) circle (4pt);
\filldraw [black] (2,3) circle (4pt);
\filldraw [black] (2,2) circle (4pt);
\filldraw [black] (1,2) circle (4pt);
\filldraw [black] (2,1) circle (4pt);
\draw[very thick] (2,6)--(2,4);
\draw[very thick] (1,4)--(2,3)--(2,1);
\end{tikzpicture}
\end{center}
An example of a CGPD for $\mathbf{r}$ is shown below. Notice its visual resemblance to $L$, with rectangles corresponding to vertices and pipes corresponding to laces.
\begin{center}
\begin{tikzpicture}[scale=0.75]
\node[loop]{
 & & & & & & & &\plaq{r}\\
 & & & & & & \plaq{r} & \plaq{h} & \\
 & & & \plaq{r} & \plaq{h} & \plaq{h} & & & \\
  & & & \plaq{j} & \plaq{} & \ecplaq{r}{red} & & & \\
  & \ecplaq{r}{violet} & \ecplaq{h}{violet} &  &  &  & & & \\
 & \ecplaq{a}{violet} & \ecplaq{h}{violet} &  &  &  & & & \\
 & \node[bgplaq,rectangle,draw,minimum size=\loopcellsize,transform shape] (plaq) {};\useasboundingbox;\begin{scope}[x=\loopcellsize,y=\loopcellsize] \draw[/linkpattern/edge,\plaqwest,\plaqnorth,color=violet] (0,0.5) .. controls (0,0.2) and (-0.2,0) .. (-0.5,0); \draw[/linkpattern/edge,\plaqeast,\plaqsouth,color=red] (0,-0.5) .. controls (0,-0.2) and (0.2,0) .. (0.5,0); \end{scope} & \ecplaq{h}{red} &  &  &  & & & \\
 \ecplaq{r}{red} & & &  &  &  & & & \\
\node[bgplaq,rectangle,draw,minimum size=\loopcellsize,transform shape] (plaq) {};\useasboundingbox;\begin{scope}[x=\loopcellsize,y=\loopcellsize] \draw[/linkpattern/edge,\plaqsouth,\plaqnorth,color=red] (0,0.5) -- (0,-0.5);
\draw[/linkpattern/edge,\plaqwest,\plaqeast,color=Cyan] (0.5,0) -- (-0.5,0); \end{scope}& &\\};
\draw (-4.5,-4.5)--(-4.5,-5.5);
\draw[blue, very thick] (2,2.5)--(1.5,2);
\draw[blue, very thick] (4,3.5)--(3.5,3);
\draw[red, very thick] (1,.5)--(-1.5,-2);
\draw[red, very thick] (-3,-2.5)--(-3.5,-3);
\draw[red, very thick] (-4,-4.5)--(-4.5,-5);
\node at (3.25,4.75) {0};
\node at (1.25,3.75) {1};
\node at (-1.75,2.75) {2};
\node at (-3.75,.75) {3};
\node at (-4.75,-2.25) {4};
\node at (-4.75,-4.25) {5};
\end{tikzpicture}
\end{center}
\end{example}
\begin{definition}
To compute the weight of a CGPD $\delta$ given a concatenated alphabet of torus weights $\mathbf{x}=\mathbf{x}^0,\mathbf{x}^1,\dots,\mathbf{x}^n$, we label the rows of rectangle $i$ top-to-bottom using the alphabet $\mathbf{x}^{i}$, and the columns left-to-right using the alphabet $\mathbf{x}^{i+1}$. Write $t_{ijk}$ for the tile with row label $x^i_j$ and column label $x^{i+1}_k$. Now the weight of $\delta$ is
\[(\mathbf{x}-\mathring{\mathbf{x}})^\delta=\prod_{i,j,k} \begin{cases}
    x^i_j-x^{i+1}_k & \textrm{if }
\begin{tikzpicture}[baseline={([yshift=-\the\dimexpr\fontdimen22\textfont2\relax]current  bounding  box.center)}]
\node[bgplaq,rectangle,draw,minimum size=\loopcellsize,transform shape] (plaq) {};\useasboundingbox;\begin{scope}[x=\loopcellsize,y=\loopcellsize] \draw[/linkpattern/edge,\plaqsouth,\plaqnorth] (0,0.5) -- (0,-0.5);
\draw[/linkpattern/edge,\plaqwest,\plaqeast,color=red] (0.5,0) -- (-0.5,0); \end{scope}
\end{tikzpicture}
, ~\plaqctr{h},~\plaqctr{v}
\textrm{ at }t_{ijk}\\\\
\hbar & \textrm{if } \begin{tikzpicture}[baseline={([yshift=-\the\dimexpr\fontdimen22\textfont2\relax]current  bounding  box.center)}]
    \node[bgplaq,rectangle,draw,minimum size=\loopcellsize,transform shape] (plaq) {};\useasboundingbox;\begin{scope}[x=\loopcellsize,y=\loopcellsize] \draw[/linkpattern/edge,\plaqwest,\plaqnorth,color=red] (0,0.5) .. controls (0,0.2) and (-0.2,0) .. (-0.5,0); \draw[/linkpattern/edge,\plaqeast,\plaqsouth] (0,-0.5) .. controls (0,-0.2) and (0.2,0) .. (0.5,0); \end{scope}\end{tikzpicture}, ~\plaqctr{j}, ~\plaqctr{r}\textrm{ at }t_{ijk}\\\\
    x^i_j - x^{i+1}_k+\hbar& \textrm{if } \plaqctr{a},~ \plaqctr{}
    \textrm{ at }t_{ijk}    
\end{cases}.\]
\end{definition}

\begin{theorem}[CGPD formula for open quiver loci]\label{thm:cgpdopen}
Given a rank array $\mathbf{r}$, we have the following expression for the CSM class of the corresponding open quiver locus:
\[\csm(\Omega_{\mathbf{r}}^\circ \subset Hom) = \sum_{\delta \in \mathcal{CGPD}(\mathbf{r})} (\mathbf{x}-\mathring{\mathbf{x}})^\delta.\]
\end{theorem}
\begin{example}
As in Example \ref{ex:a3}, let $V=(\C^1,\C^2,\C^1)$ and $\mathbf{r}=
\begin{tabular}{ccc|c}
    2 & 1 & 0 & $i \diagup j$\\\hline
  &   & 1 &0\\
  &  2 & 1 & 1\\
   1 & 1 & 1 & 2 
\end{tabular}$, so the corresponding lace array is $L=$
\begin{tikzpicture}[scale=0.5, baseline={([yshift=-.5ex]current bounding box.center)}]
\node at (1.5,2) {\small 0};
\node at (1.5,1) {\small 1};
\node at (1.5,0) {\small 2};
\filldraw [black] (1,2) circle (4pt);
\filldraw [black] (0,1) circle (4pt);
\filldraw [black] (1,1) circle (4pt);
\filldraw [black] (1,0) circle (4pt);
\draw[very thick] (1,2)--(1,1)--(1,0);
\end{tikzpicture}~. Then $\mathcal{CGPD}(L)$ contains the following three CGPDs:
\[
\begin{tikzpicture}[scale=0.75]
\node[loop]{
        &\plaq{r}&\plaq{h}\\
\plaq{r}& &\\
{
\node[bgplaq,rectangle,draw,minimum size=\loopcellsize,transform shape] (plaq) {};\useasboundingbox;\begin{scope}[x=\loopcellsize,y=\loopcellsize] \draw[/linkpattern/edge,\plaqsouth,\plaqnorth] (0,0.5) -- (0,-0.5);
\draw[/linkpattern/edge,\plaqwest,\plaqeast,color=red] (0.5,0) -- (-0.5,0); \end{scope}} &&\\};
\node at (-.75,1.1) {$a$};
\node at (1,1.75) {$b_2$};
\node at (0,1.75) {$b_1$};
\node at (-1.1,0.75) {$c$};
\node at (-1.75,0) {$b_1$};
\node at (-1.75,-1) {$b_2$};
\draw (-1.5,1.5)--(-.5,0.5);
\draw (-1.5,-1.5)--(-1.5,-2.5);
\draw[blue,very thick] (-0.5,0)--(0,0.5);
\draw[blue,very thick] (-1.5,-2)--(-1,-1.5);
\end{tikzpicture}
\hspace{.75cm}
\begin{tikzpicture}[scale=0.75]
\node[loop]{
        &\plaq{}&\plaq{r}\\
\ecplaq{h}{red}& &\\
\plaq{r}&&\\};
\node at (-.75,1.1) {$a$};
\node at (1,1.75) {$b_2$};
\node at (0,1.75) {$b_1$};
\node at (-1.1,0.75) {$c$};
\node at (-1.75,0) {$b_1$};
\node at (-1.75,-1) {$b_2$};
\draw (-1.5,1.5)--(-.5,0.5);
\draw[blue,very thick] (-0.5,-1)--(1,0.5);
\draw[blue,very thick] (-1.5,-2)--(-1,-1.5);
\draw (-1.5,-1.5)--(-1.5,-2.5);
\end{tikzpicture}\hspace{.75cm}
\begin{tikzpicture}[scale=0.75]
\node[loop]{
        &\plaq{}&\plaq{r}\\
\ecplaq{r}{red}& &\\
{
\node[bgplaq,rectangle,draw,minimum size=\loopcellsize,transform shape] (plaq) {};\useasboundingbox;\begin{scope}[x=\loopcellsize,y=\loopcellsize] \draw[/linkpattern/edge,\plaqwest,\plaqnorth,color=red] (0,0.5) .. controls (0,0.2) and (-0.2,0) .. (-0.5,0); \draw[/linkpattern/edge,\plaqeast,\plaqsouth] (0,-0.5) .. controls (0,-0.2) and (0.2,0) .. (0.5,0); \end{scope}} &&\\};
\node at (-.75,1.1) {$a$};
\node at (1,1.75) {$b_2$};
\node at (0,1.75) {$b_1$};
\node at (-1.1,0.75) {$c$};
\node at (-1.75,0) {$b_1$};
\node at (-1.75,-1) {$b_2$};
\draw (-1.5,1.5)--(-.5,0.5);
\draw (-1.5,-1.5)--(-1.5,-2.5);
\draw[blue,very thick] (-0.5,-1)--(1,0.5);
\draw[blue,very thick] (-1.5,-2)--(-1,-1.5);
\end{tikzpicture}
\]
Using Theorem \ref{thm:cgpdopen}, we compute that
\[\csm(\Omega_{\mathbf{r}}^\circ \subset Hom)= 
\hbar^2(a-b_2)(b_2-c)+\hbar^2(a-b_1+\hbar)(b_1-c)+\hbar^3(a-b_1+\hbar).\]
Note that this result agrees with Example \ref{ex:a3}.
\end{example}
\begin{proof}[Proof of Theorem \ref{thm:cgpdopen}]
This is a purely combinatorial translation of Proposition \ref{prop:pdopen}. We first note that all tiles except those on the block superantidiagonal are completely determined and may be deleted, replacing the bumps on the block diagonal with pipes continuing due southwest. There are $\ell(w_0\mathbf{w}_0)-|c(D)|$ remaining bump tiles, so we give these tiles a weight of $\hbar$. Next, note that the conditions on where pipes enter and exit the rectangles exactly encode whether the pipe dream's permutation $v$ is an element of $\perm(\mathbf{r})$. The lace array $\mathbf{s}$ is reproduced on and below the block antidiagonal of $z(\mathbf{r})$, and if the matrix of $z(\mathbf{r})$ has a 1 in block row $i$ and block column $j$ for $i \leq j$, then a pipe must enter from the right in rectangle $i$ and leave from the left in rectangle $j$.
Not requiring specific entrance or exit locations \emph{within} a particular rectangle corresponds to varying the locations of the rank jumps of $z(\mathbf{r})$ \emph{within} blocks to form $v \in \perm(\mathbf{r})$.

In fact, since we only care about what the pipes representing laces do, we can erase the other pipes; if a tile is completely blank, the erased pipes may cross or not, so its weight is the sum of a bump tile and a cross tile. If two pipes of the same color meet, whether they cross or not doesn't affect the condition from \textbf{r}, so the weights of both pipe dreams should always be included. We can therefore arbitrarily disallow same-color pipes from crossing in the diagrams while assigning the same-color bump tile the sum of a bump and a cross. This brings us to the definition of CGPDs.
\end{proof}

As a final point in this section, we consider the limit as $\hbar \to \infty$ of Theorem \ref{thm:cgpdopen}, the CGPD formula for open quiver loci. This is equivalent to setting $\hbar$ to 1 and take the lowest-degree terms. Now the only CGPDs contributing to the sum are those with the minimal number of 
$\begin{tikzpicture}
[baseline={([yshift=-\the\dimexpr\fontdimen22\textfont2\relax]current  bounding  box.center)}] \node[bgplaq,rectangle,draw,minimum size=\loopcellsize,transform shape] (plaq) {};\useasboundingbox;\begin{scope}[x=\loopcellsize,y=\loopcellsize] \draw[/linkpattern/edge,\plaqsouth,\plaqnorth] (0,0.5) -- (0,-0.5); \draw[/linkpattern/edge,\plaqwest,\plaqeast,color=red] (0.5,0) -- (-0.5,0); \end{scope} \end{tikzpicture},
\plaqctr{v},$ and $\plaqctr{h}$ tiles; for a given rank array $\mathbf{r}$, let $\mathcal{CGPD}_\infty(\mathbf{r})$ denote the set of such CGPDs.
The limiting weights of $\plaqctr{}$ and $\plaqctr{a}$ tiles correspond to interpreting these tiles as noncrossing. This is exactly the weighting from Theorem \ref{thm:newpd} on the set of reduced pipe dreams for $z(\mathbf{r})$ (which has minimal inversions in $\perm(\mathbf{r})$), so we have the following corollary of Theorem \ref{thm:cgpdopen}.
\begin{corollary}[CGPD formula for quiver polynomials]\label{cor:cgpd}
Given a CGPD $\delta$, let $c(\delta)$ denote the set of $\begin{tikzpicture}
[baseline={([yshift=-\the\dimexpr\fontdimen22\textfont2\relax]current  bounding  box.center)}] \node[bgplaq,rectangle,draw,minimum size=\loopcellsize,transform shape] (plaq) {};\useasboundingbox;\begin{scope}[x=\loopcellsize,y=\loopcellsize] \draw[/linkpattern/edge,\plaqsouth,\plaqnorth] (0,0.5) -- (0,-0.5); \draw[/linkpattern/edge,\plaqwest,\plaqeast,color=red] (0.5,0) -- (-0.5,0); \end{scope} \end{tikzpicture}, \plaqctr{v},$ and $\plaqctr{h}$ tiles in $\delta$, and recall that $t_{ijk}$ is the tile with row label $x_j^i$ and column label $x_k^{i+1}$ in $\delta$.
The quiver polynomial $[\Omega_{\mathbf{r}}]$ associated to a rank array $\mathbf{r}$ is
\[[\Omega_{\mathbf{r}}]=\sum_{\substack{\delta \in \\\mathcal{CGPD}_{\infty}(\mathbf{r})}} \prod_{t_{ijk} \in c(\delta)} 
(x^i_j-x^{i+1}_k).\]
\end{corollary}
There are the same number of terms in this sum as the sum in Theorem \ref{thm:newpd}, but it is of interest that CGPDs are simpler diagrams than pipe dreams; they look very similar to lacing diagrams, and can be enumerated directly from a lacing diagram without computing the corresponding Zelevinsky permutation.
\begin{example}
Let $V=(\C^2,\C^2,\C^1)$ and $\mathbf{r}=
\begin{tabular}{ccc|c}
    2 & 1 & 0 & $i \diagup j$\\\hline
  &   & 2 &0\\
  &  2 & 1 & 1\\
   1 & 1 & 0 & 2 
\end{tabular}$, so $\mathbf{s}=
\begin{tabular}{ccc|c}
    2 & 1 & 0 & $i \diagup j$\\\hline
  &   & 1 &0\\
  &  0 & 1 & 1\\
   0 & 1 & 0 & 2 
\end{tabular}$ and $L=$
\begin{tikzpicture}[scale=0.5, baseline={([yshift=-.5ex]current bounding box.center)}]
\node at (1.5,2) {\small 0};
\node at (1.5,1) {\small 1};
\node at (1.5,0) {\small 2};
\filldraw [black] (0,2) circle (4pt);
\filldraw [black] (1,2) circle (4pt);
\filldraw [black] (0,1) circle (4pt);
\filldraw [black] (1,1) circle (4pt);
\filldraw [black] (1,0) circle (4pt);
\draw[very thick] (1,2)--(0,1);
\draw[very thick] (1,1)--(1,0);
\end{tikzpicture}~.
We have
\[\mathcal{RP}^*(v)=\left\{
\begin{tikzpicture}[scale=0.75,baseline={([yshift=-\the\dimexpr\fontdimen22\textfont2\relax]current  bounding  box.center)}]
\node[loop]{
\gplaq{c}&\plaq{a}&\plaq{c}&\wplaq{a}&\wplaq{j}\\
\gplaq{c}&\plaq{a}&\plaq{c}&\wplaq{j}&\wplaq{}\\
\plaq{a}&\wplaq{a}&\wplaq{j}&\gplaq{}&\gplaq{}\\
\plaq{a}&\wplaq{j}&\wplaq{}&\gplaq{}&\gplaq{}\\
\wplaq{j}&\gplaq{}&\gplaq{}&\gplaq{}&\gplaq{}\\};
\node at (-2.85,2) {$a_1$};
\node at (-2.85,1) {$a_2$};
\node at (-2.85,0) {$b_1$};
\node at (-2.85,-1) {$b_2$};
\node at (-2.85,-2) {$c$};
\node at (-2,2.85) {$c$};
\node at (-1,2.85) {$b_1$};
\node at (0,2.85) {$b_2$};
\node at (1,2.85) {$a_1$};
\node at (2,2.85) {$a_2$};
\end{tikzpicture},
\hspace{.25cm}
\begin{tikzpicture}[scale=0.75,baseline={([yshift=-\the\dimexpr\fontdimen22\textfont2\relax]current  bounding  box.center)}]
\node[loop]{
\gplaq{c}&\plaq{a}&\plaq{a}&\wplaq{a}&\wplaq{j}\\
\gplaq{c}&\plaq{c}&\plaq{a}&\wplaq{j}&\wplaq{}\\
\plaq{a}&\wplaq{a}&\wplaq{j}&\gplaq{}&\gplaq{}\\
\plaq{c}&\wplaq{j}&\wplaq{}&\gplaq{}&\gplaq{}\\
\wplaq{j}&\gplaq{}&\gplaq{}&\gplaq{}&\gplaq{}\\};
\node at (-2.85,2) {$a_1$};
\node at (-2.85,1) {$a_2$};
\node at (-2.85,0) {$b_1$};
\node at (-2.85,-1) {$b_2$};
\node at (-2.85,-2) {$c$};
\node at (-2,2.85) {$c$};
\node at (-1,2.85) {$b_1$};
\node at (0,2.85) {$b_2$};
\node at (1,2.85) {$a_1$};
\node at (2,2.85) {$a_2$};
\end{tikzpicture},
\hspace{.25cm}
\begin{tikzpicture}[scale=0.75,baseline={([yshift=-\the\dimexpr\fontdimen22\textfont2\relax]current  bounding  box.center)}]
\node[loop]{
\gplaq{c}&\plaq{a}&\plaq{c}&\wplaq{a}&\wplaq{j}\\
\gplaq{c}&\plaq{a}&\plaq{a}&\wplaq{j}&\wplaq{}\\
\plaq{a}&\wplaq{a}&\wplaq{j}&\gplaq{}&\gplaq{}\\
\plaq{c}&\wplaq{j}&\wplaq{}&\gplaq{}&\gplaq{}\\
\wplaq{j}&\gplaq{}&\gplaq{}&\gplaq{}&\gplaq{}\\};
\node at (-2.85,2) {$a_1$};
\node at (-2.85,1) {$a_2$};
\node at (-2.85,0) {$b_1$};
\node at (-2.85,-1) {$b_2$};
\node at (-2.85,-2) {$c$};
\node at (-2,2.85) {$c$};
\node at (-1,2.85) {$b_1$};
\node at (0,2.85) {$b_2$};
\node at (1,2.85) {$a_1$};
\node at (2,2.85) {$a_2$};
\end{tikzpicture}
\right\}\]
while
\[\mathcal{CGPD}_{\infty}(v)=\left\{
\begin{tikzpicture}[scale=0.75,baseline={([yshift=-\the\dimexpr\fontdimen22\textfont2\relax]current  bounding  box.center)}]
\node[loop]{
&\ecplaq{r}{red}&\ecplaq{h}{red}\\
&\ecplaq{a}{red}&\ecplaq{h}{red}\\
\ecplaq{r}{red}&&\\
\node[bgplaq,rectangle,draw,minimum size=\loopcellsize,transform shape] (plaq) {};\useasboundingbox;\begin{scope}[x=\loopcellsize,y=\loopcellsize] \draw[/linkpattern/edge,\plaqwest,\plaqnorth,color=red] (0,0.5) .. controls (0,0.2) and (-0.2,0) .. (-0.5,0); \draw[/linkpattern/edge,\plaqeast,\plaqsouth] (0,-0.5) .. controls (0,-0.2) and (0.2,0) .. (0.5,0); \end{scope}&&\\};
\node at (-.75,1.5) {$a_1$};
\node at (-.75,.6) {$a_2$};
\node at (-1.1,.25) {$c$};
\node at (0,2.25) {$b_1$};
\node at (1,2.25) {$b_2$};
\node at (-1.75,-.5) {$b_1$};
\node at (-1.75,-1.5) {$b_2$};
\draw (-1.5,-2)--(-1.5,-3);
\draw (-.5,0)--(-1.5,1);
\draw[color=red,very thick] (0,0) -- (-.5,-.5);
\draw[color=blue,very thick] (-1,-2) -- (-1.5,-2.5);
\end{tikzpicture},
\hspace{1cm}
\begin{tikzpicture}[scale=0.75,baseline={([yshift=-\the\dimexpr\fontdimen22\textfont2\relax]current  bounding  box.center)}]
\node[loop]{
&\plaq{}&\ecplaq{r}{red}\\
&\ecplaq{h}{red}&\ecplaq{a}{red}\\
\plaq{r}&&\\
\node[bgplaq,rectangle,draw,minimum size=\loopcellsize,transform shape] (plaq) {};\useasboundingbox;\begin{scope}[x=\loopcellsize,y=\loopcellsize] \draw[/linkpattern/edge,\plaqsouth,\plaqnorth] (0,0.5) -- (0,-0.5); \draw[/linkpattern/edge,\plaqwest,\plaqeast,color=red] (0.5,0) -- (-0.5,0); \end{scope}&&\\};
\node at (-.75,1.5) {$a_1$};
\node at (-.75,.6) {$a_2$};
\node at (-1.1,.25) {$c$};
\node at (0,2.25) {$b_1$};
\node at (1,2.25) {$b_2$};
\node at (-1.75,-.5) {$b_1$};
\node at (-1.75,-1.5) {$b_2$};
\draw (-1.5,-2)--(-1.5,-3);
\draw (-.5,0)--(-1.5,1);
\draw[color=red,very thick] (1,0) -- (-.5,-1.5);
\draw[color=blue,very thick] (-1,-2) -- (-1.5,-2.5);
\end{tikzpicture},
\hspace{1cm}
\begin{tikzpicture}[scale=0.75,baseline={([yshift=-\the\dimexpr\fontdimen22\textfont2\relax]current  bounding  box.center)}]
\node[loop]{
&\ecplaq{r}{red}&\ecplaq{h}{red}\\
&\ecplaq{j}{red}&\ecplaq{r}{red}\\
\plaq{r}&&\\
\node[bgplaq,rectangle,draw,minimum size=\loopcellsize,transform shape] (plaq) {};\useasboundingbox;\begin{scope}[x=\loopcellsize,y=\loopcellsize] \draw[/linkpattern/edge,\plaqsouth,\plaqnorth] (0,0.5) -- (0,-0.5); \draw[/linkpattern/edge,\plaqwest,\plaqeast,color=red] (0.5,0) -- (-0.5,0); \end{scope}&&\\};
\node at (-.75,1.5) {$a_1$};
\node at (-.75,.6) {$a_2$};
\node at (-1.1,.25) {$c$};
\node at (0,2.25) {$b_1$};
\node at (1,2.25) {$b_2$};
\node at (-1.75,-.5) {$b_1$};
\node at (-1.75,-1.5) {$b_2$};
\draw (-1.5,-2)--(-1.5,-3);
\draw (-.5,0)--(-1.5,1);
\draw[color=red,very thick] (1,0) -- (-.5,-1.5);
\draw[color=blue,very thick] (-1,-2) -- (-1.5,-2.5);
\end{tikzpicture}
\right\}.\]
Both Theorem \ref{thm:newpd} and Corollary \ref{cor:cgpd} give
\[[\Omega_{\mathbf{r}}]=(a_1-b_2)(a_2-b_2)+(a_2-b_1)(b_2-c)+(a_1-b_2)(b_2-c).\]
\end{example}
\FloatBarrier

\printbibliography
\end{document}